\documentclass[a4paper,12pt]{article}
\usepackage{graphicx}
\usepackage{psfrag}

\usepackage{mathrsfs}
\usepackage{amssymb}
\usepackage{amsfonts}
\usepackage{latexsym}
\usepackage{amsmath}
\usepackage{amsthm}
\usepackage[cp1251]{inputenc}
\usepackage[english]{babel}

\newcommand{\er}[1]{\textrm{(\ref{#1})}}
\def\lb{\label}

\theoremstyle{plain}
\newtheorem{Def}{\bf Definition}
\newtheorem{theorem}{\bf Theorem}[section]
\newtheorem{lemma}[theorem]{\bf Lemma}
\newtheorem{proposition}[theorem]{\bf Proposition}

\theoremstyle{remark}

\renewcommand{\a}{\alpha}

\newcommand{\g}{\gamma}           

\newcommand{\G}{\Gamma}           

\renewcommand{\d}{\delta}           

\newcommand{\D}{\Delta}           

\newcommand{\ve}{\varepsilon}     

            \newcommand{\cH}{\mathcal{H}}

\newcommand{\e}{\eta}             

\newcommand{\vt}{\vartheta}

\renewcommand{\l}{\lambda}          \newcommand{\cM}{\mathcal{M}}

\newcommand{\m}{\mu}              

\newcommand{\n}{\nu}              

\renewcommand{\r}{\rho}             

\newcommand{\s}{\sigma}

\renewcommand{\t}{\tau}

\newcommand{\vp}{\varphi}

\newcommand{\p}{\psi}             
 
             \newcommand{\cZ}{\mathcal{Z}}
      
\renewcommand{\o}{\omega}
\renewcommand{\O}{\Omega}
\newcommand{\x}{\xi}

\newcommand{\vk}{\varkappa}

  \def\mH{{\mathscr H}}

  \def\mV{{\mathscr V}}

\newcommand{\gD}{\mathfrak{D}}

\newcommand{\gR}{\mathfrak{R}}
\newcommand{\gS}{\mathfrak{S}}

\def\Z{\mathbb{Z}}
\def\R{\mathbb{R}}
\def\C{\mathbb{C}}

\def\N{\mathbb{N}}

\def\qqq{\qquad}
\def\qq{\quad}

\let\ge\geqslant
\let\le\leqslant
\let\geq\geqslant

\newcommand{\ca}{\begin{cases}}
\newcommand{\ac}{\end{cases}}
\newcommand{\ma}{\begin{pmatrix}}
\newcommand{\am}{\end{pmatrix}}

\def\lt{\biggl}
\def\rt{\biggr}
\let\geq\geqslant

\renewcommand{\[}{\begin{equation}}
\renewcommand{\]}{\end{equation}}

\def\wt{\widetilde}

\def\sm{\setminus}
\def\es{\emptyset}
\def\no{\noindent}
\def\ol{\overline}
\def\iy{\infty}

\def\/{\over}
\def\ts{\times}

\def\os{\oplus}
\def\ss{\subset}

\def\Im{\mathop{\rm Im}\nolimits}

\def\sign{\mathop{\rm sign}\nolimits}

\def\const{\mathop{\rm const}\nolimits}

\def\BBox{\hspace{1mm}\vrule height6pt width5.5pt depth0pt \hspace{6pt}}

\begin{document}
\title{Schr\"odinger operators on
armchair nanotubes. II}
\author{
Andrey Badanin
\begin{footnote} {
Arkhangelsk State Technical University,
e-mail: a.badanin@agtu.ru }
\end{footnote}
 \and
 Jochen Br\"uning
 \begin{footnote} {
 Institut f\"ur Mathematik, Humboldt Universit\"at zu Berlin,
 e-mail: bruening@math.hu-berlin.de }
 \end{footnote}
\and
Evgeny Korotyaev
\begin{footnote} {
Institut f\"ur  Mathematik,  Humboldt Universit\"at zu Berlin,
e-mail: evgeny@math.hu-berlin.de\ \
}
\end{footnote}
}

\maketitle

\begin{abstract}
\no We consider the Schr\"odinger operator with a periodic potential on
quasi-1D models of armchair single-wall nanotubes.
The spectrum of this operator consists of an absolutely continuous part
(intervals separated by gaps) plus an infinite number of eigenvalues
with infinite multiplicity.
We describe the absolutely continuous spectrum of
the Schr\"odinger operator: 1) the multiplicity,
2) endpoints of the gaps, they  are given by
periodic or antiperiodic eigenvalues or resonances (branch points
of the Lyapunov function),
3) resonance gaps, where the Lyapunov function is non-real.
We determine the asymptotics of the gaps at high energy.
\end{abstract}

\section{Introduction  and main results}
\setcounter{equation}{0}

\begin{figure}
\centering
\noindent
\tiny
\hfill
\includegraphics[height=.4\textheight]{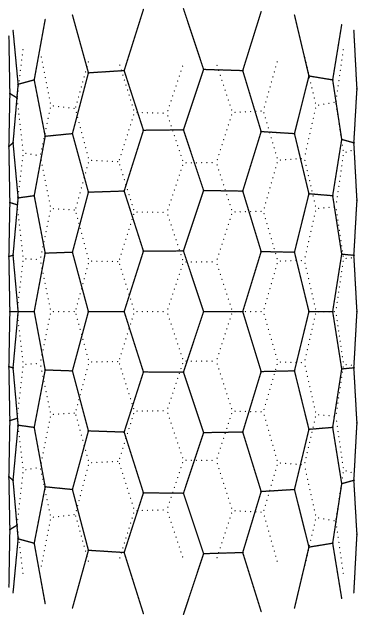}
\hfill {
\psfrag{g010}{$\Gamma_{0,1}$}
\psfrag{g020}{$\Gamma_{0,2}$}
\psfrag{g030}{$\Gamma_{0,3}$}
\psfrag{g040}{$\Gamma_{0,4}$}
\psfrag{g050}{$\Gamma_{0,5}$}
\psfrag{g060}{$\Gamma_{0,6}$}
\psfrag{g110}{$\Gamma_{1,1}$}
\psfrag{g120}{$\Gamma_{1,2}$}
\psfrag{g130}{$\Gamma_{1,3}$}
\psfrag{g140}{$\Gamma_{1,4}$}
\psfrag{g150}{$\Gamma_{1,5}$}
\psfrag{g160}{$\Gamma_{1,6}$}
\psfrag{g-110}{$\Gamma_{-1,1}$}
\psfrag{g-120}{$\Gamma_{-1,2}$}
\psfrag{g-130}{$\Gamma_{-1,3}$}
\psfrag{g-140}{$\Gamma_{-1,4}$}
\psfrag{g-150}{$\Gamma_{-1,5}$}
\psfrag{g-160}{$\Gamma_{-1,6}$}
\includegraphics[height=.4\textheight]{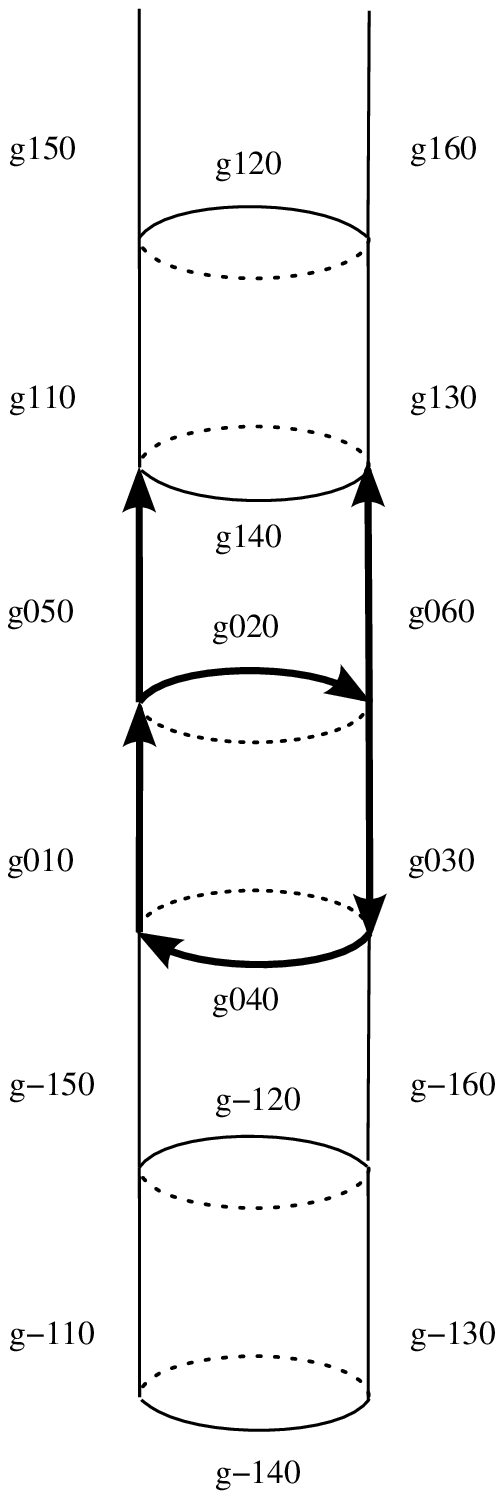}
}
\hfill~
\caption{Armchair graph for $N=10$ and for $N=1$.}
\lb{fig1}
\end{figure}

Consider the Schr\"odinger operator $\mH =-\D+\mV_q$
 with a periodic potential $\mV_q$ on so called armchair
graph $\G^{N},N\ge 1$.
In order to describe the graph $\G^N$, we introduce the fundamental
cell $\wt\G=\cup_1^6\wt\G_j \ss\R^2$, where
$\wt\G_j=\{x=\wt{\bf r}_j+t{\bf e}_j, t\in [0,1]\},j\in\N_6$ is the edge of length 1, and
\begin{multline}
\N_m=\{1,2,..,m\},\qq
{\bf e}_1={\bf e}_6={1\/2}(1,\sqrt 3),\qq
{\bf e}_2={\bf e}_4=(1,0),\qq {\bf e}_3=-{\bf e}_5={1\/2}(1,-\sqrt 3),\\
\wt{\bf r}_1=(0,0),\qq \wt{\bf r}_2=\wt{\bf r}_5=\wt{\bf r}_1+{\bf e}_1,\qq
\wt{\bf r}_3=\wt{\bf r}_6=\wt{\bf r}_2+{\bf e}_2,\qq
\wt{\bf r}_4=\wt{\bf r}_3+{\bf e}_3.
\end{multline}
We define the strip graph $\wt\G^N$ by
$$
\wt\G^N=\cup_{(n,k)\in\Z\ts\N_N}
(\wt\G+k{\bf e}_h+n{\bf e}_v)\ss\R^2,
\qq {\bf e}_h=(3,0),\qq {\bf e}_v=(0,\sqrt 3).
$$
Vertices of $\wt\G^N$ are
$\wt {\bf r}_j+k{\bf e}_h+n{\bf e}_v,(n,j,k)\in\Z\ts\N_6\ts\N_N$.
If we identify the vertices
$\wt{\bf r}_1+n{\bf e}_v$ and $\wt{\bf r}_1+N{\bf e}_h+n{\bf e}_v$
of $\wt\G^N$
for each $n\in\Z$, then we obtain the graph $\G^N$, given by
$$
\G^N=\cup_{\o\in \cZ} \G_\o,\qq
\o=(n,j,k)\in \cZ=\Z\ts \N_6\ts \Z_N,
\qq \Z_N=\Z/(N\Z),
$$
where $\G_{\o}=\wt\G_j+k{\bf e}_h+n{\bf e}_v$,
see Fig. \ref{fig1}, \ref{fig2}.
Let ${\bf r}_\o=\wt{\bf r}_j+k{\bf e}_h+n{\bf e}_v$
be a starting point of the edge $\G_{\o}$.
We have the coordinate
$x={\bf r}_\o+t{\bf e}_j$
and the local coordinate $t\in [0,1]$ on $\G_\o$.
Thus we give an orientation on the edge.
For each function $y$ on $\G^N$ we define a function $y_\o=y|_{\G_\o}, \o\in \cZ$. We identify each function $y_\o$ on $\G_\o$ with a function on $[0,1]$ by using the local coordinate $t\in [0,1]$. Define the Hilbert space $L^2(\G^N)=\os_{\o\in \cZ} L^2(\G_\o)$.
Let $C(\G^N)$ be the space   of continuous functions on $\G^N$.
We define the Sobolev space $W^2(\G^N)$ that consists of all
functions $y=(y_\o)_{\o\in\cZ}\in L^2(\G^N),
(y_\o'')_{\o\in\cZ}\in L^2(\G^N)$ and satisfy

\no {\bf   Kirchhoff Boundary Conditions:} {\it
$y\in C(\G^N)$ satisfies for each vertex $A$ of $\G^N$

\[
\lb{KirC}
\sum_{\o\in E_A}(-1)^b y_\o'(b)=0,\qq \text{where}\qq
E_A=\{\o\in\cZ:A\in\G_\o\},\qq b=b(\o,A),
\]
where if $A={\bf r}_\o$ is a starting point of $\G_\o$
$($i.e. $t=0$ at $A)$,
then $b(\o,A)=0$,

\no \qqq \ \ if $A={\bf r}_\o+{\bf e_j}$ is an endpoint of $\G_\o$
$($i.e. $t=1$ at $A)$,
then $b(\o,A)=1$.
}

The Kirchhoff Conditions \er{KirC} mean  that
the sum of derivatives of $y$ at each vertex of $\G^N$ equals 0 and the orientation of edges gives
the sign $\pm$. Our operator $\mH$ on $\G^N$ acts in the Hilbert space
$L^2(\G^N)$ and is given  by
$(\mH y)_\o=-y_\o''+q y_\o$, where $\ y=(y_\o)_{\o\in \cZ}\in
\gD(\mH)=W^2(\G^N)$  and $(\mV_q y)_\o=qy_\o, q\in L^2(0,1)$.
If the potential $q$ is {\bf even}, i.e., $q\in L^2_{even}(0,1)=\{q\in
L^2(0,1): q(t)=q(1-t),t\in[0,1]\}$, then
the orientation of edges is not important.
The standard arguments (see \cite{KL}) yield  that $\mH$ is
self-adjoint.

The considered model was introduced by Pauling \cite{Pa}
and was systematically developed in the series of
articles by Ruedenberg and Scherr \cite{RS}.
Further progress is discussed in \cite{KL},\cite{KL1},\cite{BBKL}, \cite{Ha}, \cite{SDD}  and see references therein.

\begin{figure}
\centering
\noindent
{
\tiny
\psfrag{g011}{$\Gamma_{0,1,1}$}
\psfrag{g021}{$\Gamma_{0,2,1}$}
\psfrag{g031}{$\Gamma_{0,3,1}$}
\psfrag{g041}{$\Gamma_{0,4,1}$}
\psfrag{g051}{$\Gamma_{0,5,1}$}
\psfrag{g061}{$\Gamma_{0,6,1}$}
\psfrag{g012}{$\Gamma_{0,1,2}$}
\psfrag{g022}{$\Gamma_{0,2,2}$}
\psfrag{g032}{$\Gamma_{0,3,2}$}
\psfrag{g042}{$\Gamma_{0,4,2}$}
\psfrag{g052}{$\Gamma_{0,5,2}$}
\psfrag{g062}{$\Gamma_{0,6,2}$}
\psfrag{g01N}{$\Gamma_{0,1,N}$}
\psfrag{g02N}{$\Gamma_{0,2,N}$}
\psfrag{g03N}{$\Gamma_{0,3,N}$}
\psfrag{g04N}{$\Gamma_{0,4,N}$}
\psfrag{g05N}{$\Gamma_{0,5,N}$}
\psfrag{g06N}{$\Gamma_{0,6,N}$}
\psfrag{g111}{$\Gamma_{1,1,1}$}
\psfrag{g121}{$\Gamma_{1,2,1}$}
\psfrag{g131}{$\Gamma_{1,3,1}$}
\psfrag{g141}{$\Gamma_{1,4,1}$}
\psfrag{g151}{$\Gamma_{1,5,1}$}
\psfrag{g161}{$\Gamma_{1,6,1}$}
\psfrag{g112}{$\Gamma_{1,1,2}$}
\psfrag{g122}{$\Gamma_{1,2,2}$}
\psfrag{g132}{$\Gamma_{1,3,2}$}
\psfrag{g142}{$\Gamma_{1,4,2}$}
\psfrag{g152}{$\Gamma_{1,5,2}$}
\psfrag{g162}{$\Gamma_{1,6,2}$}
\psfrag{g11N}{$\Gamma_{1,1,N}$}
\psfrag{g2N}{$\Gamma_{1,2,N}$}
\psfrag{g13N}{$\Gamma_{1,3,N}$}
\psfrag{g14N}{$\Gamma_{1,4,N}$}
\psfrag{g15N}{$\Gamma_{1,5,N}$}
\psfrag{g16N}{$\Gamma_{1,6,N}$}
\psfrag{g-111}{$\Gamma_{-1,1,1}$}
\psfrag{g-121}{$\Gamma_{-1,2,1}$}
\psfrag{g-131}{$\Gamma_{-1,3,1}$}
\psfrag{g-141}{$\Gamma_{-1,4,1}$}
\psfrag{g-151}{$\Gamma_{-1,5,1}$}
\psfrag{g-161}{$\Gamma_{-1,6,1}$}
\psfrag{g-112}{$\Gamma_{-1,1,2}$}
\psfrag{g-122}{$\Gamma_{-1,2,2}$}
\psfrag{g-132}{$\Gamma_{-1,3,2}$}
\psfrag{g-142}{$\Gamma_{-1,4,2}$}
\psfrag{g-152}{$\Gamma_{-1,5,2}$}
\psfrag{g-162}{$\Gamma_{-1,6,2}$}
\psfrag{g-11N}{$\Gamma_{-1,1,N}$}
\psfrag{g-12N}{$\Gamma_{-1,2,N}$}
\psfrag{g-13N}{$\Gamma_{-1,3,N}$}
\psfrag{g-14N}{$\Gamma_{-1,4,N}$}
\psfrag{g-15N}{$\Gamma_{-1,5,N}$}
\psfrag{g-16N}{$\Gamma_{-1,6,N}$}
\psfrag{O1}{$\Omega_{1}$}
\psfrag{O2}{$\Omega_{2}$}
\includegraphics[width=.8\textwidth,height=.5\textwidth]{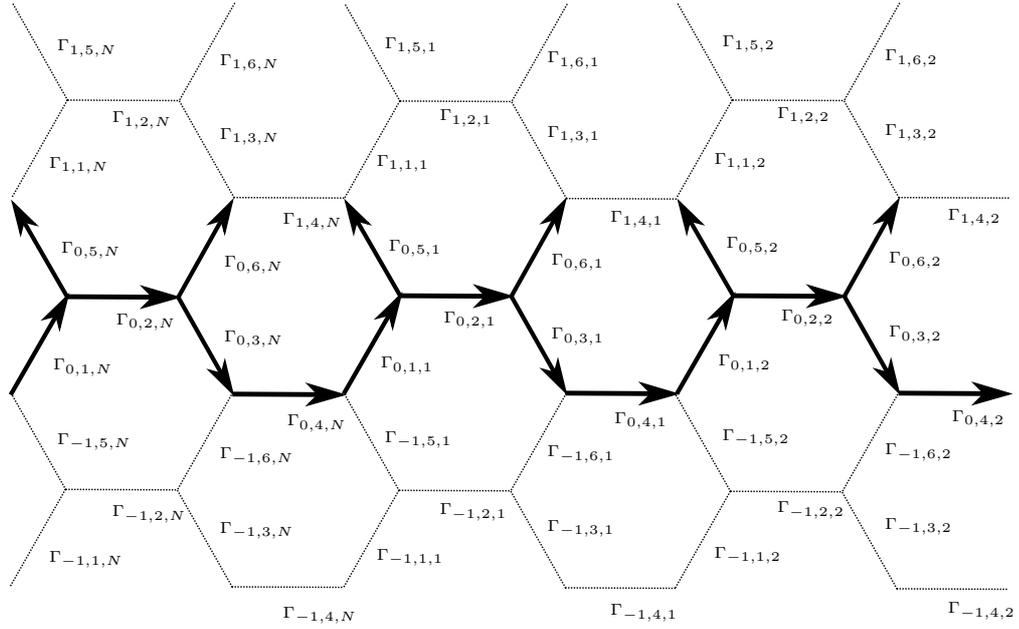}
}
\caption{ A piece of a nanotube $\G^N$.
The fundamental domain is marked by a bold line.
}
\lb{fig2}
\end{figure}

\begin{figure}
\noindent
\centering
\tiny
\psfrag{A}[r][r]{$A_1$}
\psfrag{B}[r][r]{$A_2$}
\psfrag{a}{$\O_1$}
\psfrag{b}{$\O_2$}
\psfrag{c}{$\O$}
\psfrag{(1,0)}{$(1,0)$}
\psfrag{(2,0)}{$(2,0)$}
\psfrag{(3,0)}{$(3,0)$}
\psfrag{(4,0)}{$(4,0)$}
\psfrag{(1,1)}{$(1,1)$}
\psfrag{(2,1)}{$(2,1)$}
\psfrag{(3,1)}{$(3,1)$}
\psfrag{(2,2)}{$(2,2)$}
\psfrag{(3,2)}{$(3,2)$}
\includegraphics[width=.5\textwidth]{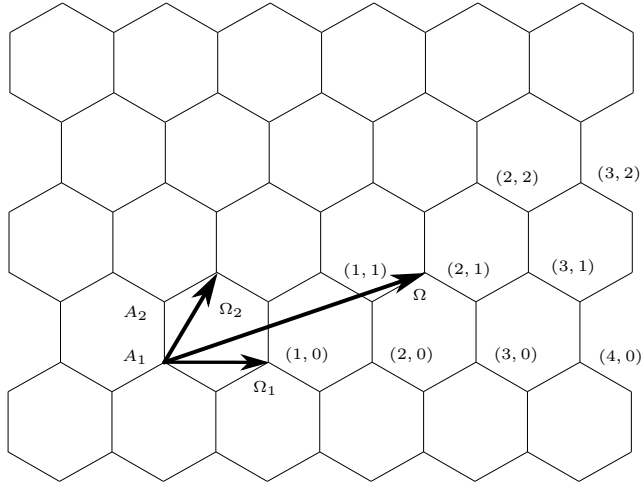}
\caption{The honeycomb lattice of nanotube.
The unit cell is spanned by the vectors $\O_1$ and $\O_2$.
The type of the nanotube is defined by the pair $(N_1,N_2)\in\N^2$,
$N_1\geq N_2$, and corresponding chiral vector
$\O=N_1\O_1+N_2\O_2$.
}
\label{fig5}
\end{figure}

For  the convenience of the reader we  briefly describe the structure
of carbon
nanotubes, see \cite{Ha}, \cite{SDD}. Graphene is a single 2D layer
of graphite
forming a honeycomb lattice, see Fig. \ref{fig5}.
A carbon nanotube is a honeycomb lattice "rolled up" into a
cylinder, see Fig. \ref{fig1}.
In carbon nanotubes, the graphene sheet is "rolled up" in such a way
that the so-called chiral  vector $\O=N_1\O_1+N_2\O_2$ becomes
the circumference of the tube, where $\O_1, \O_2$ are defined in
 Fig \ref{fig5}. The chiral  vector $\O$, which is usually denoted by the pair
of integers $(N_1,N_2)$, uniquely defines a particular tube.
Tubes of type $(N,0)$ are called zigzag tubes.
$(N,N)$-tubes are called armchair tubes.

Recall the needed properties of the Hill operator $\wt H
y=-y''+q(t)y$ on the real line with a periodic potential
$q(t+1)=q(t),t\in \R$. The spectrum of $\wt H$ is purely absolutely
continuous and consists of intervals
$\wt\s_n=[\wt\l_{n-1}^+,\wt\l_n^-], n\ge 1$. These intervals are
separated by the gaps $\wt\g_n=(\wt\l_n^-,\wt\l_n^+)$ of length
$|\wt\g_n|\ge 0$. If a gap $\wt\g_n$ is degenerate, i.e.
$|\wt\g_n|=0$, then the corresponding segments $\wt\s_n,\wt\s_{n+1}$
merge. For the equation $-y''+q(t)y=\l y$ on the real line we define
the fundamental solutions $\vt(t,\l)$ and $\vp(t,\l),t\in \R$
satisfying $\vt(0,\l)=\vp'(0,\l)=1, \vt'(0,\l)=\vp(0,\l)=0$.
We define the monodromy matrix $\wt\cM$, the Lyapunov function $F$,
and the function $F_-$ by
\[
\lb{Ho}
\wt\cM=\ma\vt_1 & \vp_1 \\
\vt_1' & \vp_1'\am,\qq
F={\vp_1'+\vt_1\/2},\qq  F_-={\vp_1'-\vt_1\/2},
\]
where
$\vp_1=\vp(1,\cdot),\vt_1=\vt(1,\cdot),\vp_1'=\vp'(1,\cdot),
\vt_1'=\vt'(1,\cdot)$.
The function $F$ has only simple zeros $\e_n,n\ge 1$, which satisfy
$\e_1<\e_2<...$
The sequence $\wt\l_0^+<\wt\l_1^-\le \wt\l_1^+\ <...$ is the
spectrum of the equation $-y''+qy=\l y$ with 2-periodic boundary
conditions,  that is  $y(t+2)=y(t), t\in \R$.
 Here equality $\wt\l_n^-= \wt\l_n^+$ means that $\wt\l_n^\pm$ is an eigenvalue of
multiplicity 2. Note that $F(\wt\l_{n}^{\pm})=(-1)^n, \  n\ge 1$.
The lowest  eigenvalue $\wt\l_0^+$ is simple, $F(\wt\l_0^+)=1$, and
the corresponding eigenfunction has period 1. The eigenfunctions
corresponding to $\wt\l_n^{\pm}$ have period 1 if $n$ is even, and
they are anti-periodic, that is $y(t+1)=-y(t),\ t\in \R$, if $n$ is
odd. The derivative of the Lyapunov function has a zero $\wt\l_n$ in
each interval $[\l^-_n,\l^+_n]$, that is $ F'(\wt\l_n)=0$.  Let
$\m_n, n\ge 1,$ be the spectrum of the problem $-y''+qy=\l y,
y(0)=y(1)=0$ (the Dirichlet spectrum).
Define the set $\s_D=\{\m_n, n\ge 1\}$ and note
that $\s_D=\{\l\in\C: \vp(1,\l)=0\}$. It is well-known that $\m_n
\in [\wt\l^-_n,\wt\l^+_n ],n\ge 1$.

For simplicity we shall denote $\G_{\a,1}\ss \G^1$
by  $\G_{\a}$,
for $\a=(n,j)\in \cZ_1=\Z\ts \N_6$.
Thus $\G^1=\cup_{\a\in \cZ_1}
\G_\a$, see Fig \ref{fig1}.
We introduce the self adjoint operator $H_k$ acting in the Hilbert space
 $L^2(\G^1)$ and given by
$(H_k f)_\a=-f_\a''+q f_\a,\ (f_\a)_{\a\in \cZ_1},
(f_\a'')_{\a\in
\cZ_1}\in L^2(\G^1)$, where the components
$f_\a,\a\in \cZ_1$  satisfy
{\bf the Kirchhoff conditions}:
\begin{multline}
\label{1K0} f_{n,1}(1)=f_{n,2}(0)=f_{n,5}(0),\qqq
f_{n,2}(1)=f_{n,3}(0)=f_{n,6}(0),\\
f_{n,3}(1)=f_{n,4}(0)=f_{n-1,6}(1),\qq s^k
f_{n,4}(1)=f_{n,1}(0)=f_{n-1,5}(1),\qq s=e^{i{2\pi \/N}},
\end{multline}
\begin{multline}
\label{1K1}
f_{n,1}'(1)-f_{n,2}'(0)-f_{n,5}'(0)=0,\qq
f_{n,2}'(1)-f_{n,3}'(0)-f_{n,6}'(0)=0,\\
f_{n,3}'(1)-f_{n,4}'(0)+f_{n-1,6}'(1)=0,\qq
s^kf_{n,4}'(1)-f_{n,1}'(0)+f_{n-1,5}'(1)=0.
\end{multline}

The operator $H_k$ has four Floquet solutions
$\p_{k}^{\n,\pm}=(\p_{k,\a}^{\n,\pm})_{\a\in\cZ_1},\n=1,2$ satisfying
the condition
$
\ma
\p_{k,1,5}^{\n,\pm}(1)\\ \p_{k,1,6}^{\n,\pm}(1)\am
=\t_{k,\n}^{\pm 1}
\ma \p_{k,0,5}^{\n,\pm}(1)\\ \p_{k,0,6}^{\n,\pm}(1)\am.
$
For each $k\in\Z_N$ we introduce two Lyapunov functions
$F_{k,\n}={1\/2}(\t_{k,\n}+\t_{k,\n}^{-1}),\n\in\N_2$.
Recall the results from \cite{BBKL}:

\no {\it
 The operator $\mH$ is unitarily equivalent to
$H=\os_1^N H_k$.
The  following identities hold true:
\begin{multline}
\lb{T3-1}
\s(H_k)=\s_{\iy}(H_k)\cup\s_{ac}(H_k),\qq  \s_{\iy}(H_k)=\s_D,\\
\s_{ac}(H_k)=\{\l\in\R: F_{k,\nu}(\l)\in [-1,1]\ \text{for some}\
\nu\in\N_2\},
\end{multline}
\[
\lb{DeLk}
F_{k,\nu}=\xi_k-(-1)^\n\sqrt{\r_k},\ \n=1,2,\qq
\xi_k={9F^2-F_-^2-1\/2}-s_k^2,\qq
\r_k=(9F^2-s_k^2)c_k^2+s_k^2F_-^2
\]
for each $k\in \Z_N$, where
$
s_k=\sin {\pi k\/N},\qq
c_k=\cos {\pi k\/N}.
$
Here the functions $F_{k,1},F_{k,2}$ are branches
of the Lyapunov functions $F_k=\xi_k+\sqrt{\r_k}$, analytic
on the two sheeted Riemann surface
$\gR_k$ defined by $\sqrt {\r_k}$.
}

\no {\bf Remark.}
We take the branch of $\sqrt{\r_k}$ such that
$\sqrt{\r_k(\l)}>0$, where $\r_k(\l)> 0,\l\in\R$.
Then $F_{k,1}=\xi_k+\sqrt{\r_k}>
F_{k,2}=\xi_k-\sqrt{\r_k}$ for such $\l$.
Note that $F_{k,\n},k\not\in\{0,{N\/2}\}$
have branch points on the real line.
The functions $F_{0,\n}',\n=1,2$
have steps at the points $\e_n,n\ge 1$. The functions
$F_{m,\n}',m={N\/2}\in\Z$  have steps at the zeros of $F_-$.
Note that using other branches of $\sqrt{\r_k}$
we could obtain a new smooth functions $F_{k,\n}$ on real axis
for $k\in\{0,{N\/2}\}$, but this choice is not
convenient for our proof.

 We define the entire functions
\[
\lb{S3b}
D_k^\pm=4(F_{k,1}\mp 1)(F_{k,2}\mp 1).
\]
\no The zeros $\l_{\n,2n}^{k,\pm}, n\ge 0,\n=1,2$, of the function
$D_k^+$ are the periodic eigenvalues.
The zeros $\l_{\n,2n-1}^{k,\pm},n\ge 1,\n=1,2$, of
$D_k^-$ are the antiperiodic eigenvalues.
Let $\l_{2,0}^{k,+}\le\l_{1,0}^{k,+}
\le\l_{1,2}^{k,-}\le\l_{2,2}^{k,-}\le
\l_{2,2}^{k,+}\le\l_{1,2}^{k,+}\le...$
and
$
\l_{2,1}^{k,-}\le\l_{1,1}^{k,-}\le\l_{1,1}^{k,+}\le\l_{2,1}^{k,+}
\le\l_{2,3}^{k,-}\le\l_{1,3}^{k,-}\le...
$
counted with multiplicities.
This labeling is convenient for us and associated with
the Lyapunov functions $F_{k,1},F_{k,2}$ (see Fig.\ref{lfd}).

A zero of $\r_k, k\in \Z_N$ is called a {\bf resonance} of $H_k$.
Roughly speaking the simple real resonances create gaps.
There exist real and non-real resonances for $k\not\in\{0,{N\/2}\}$
(see \cite{BBKL}). Note that in the case of
zigzag nanotube all resonances are real
 \cite{KL}, \cite{KL1}.

We define the functions
\[
\lb{uv}
u_k=|F_-|-s_k^2,\qq v_k=|F_-|-c_k^2,\qq \qq k\in\Z_N.
\]

\begin{theorem}
\lb{Tk}  Let $k\in \Z_N$. Then
the identity $\s_{ac}(H_k)=\cup_{\n\in\N_2,n\ge 1}S_{\n,n}^k$ holds,
where
the spectral bands $S_{\n,n}^k=[E_{\n,n-1}^{k,+},E_{\n,n}^{k,-}],n\ge 1,
\n=1,2$
satisfy:
\[
\lb{s}
E_{\n,p-1}^{k,\pm}=\l_{\n,p-1}^{k,\pm},\qq
E_{2,p}^{k,\pm}=\l_{2,p}^{0,\pm},\qq
E_{1,p}^{k,\pm}=\ca
\l_{1,p}^{0,\pm}\qq \text{if}\qq
v_k(\l_{1,p}^{0,\pm})\ge 0\\
r_{k,n}^\pm\qq\ \text{if}\qq
v_k(\l_{1,p}^{0,\pm})<0
\ac\!\!\!,\qq p=2n-1,
\]
$E_{1,p}^{k,\pm}=r_{k,n}^\pm$ for $k\not\in\{0,{N\/2}\}$
and for large $n\ge 1$, where $r_{k,n}^\pm$ are given by
$$
r_{k,n}^-=\min\{\l\in\ol\vk_n:\r_k(\l)=0\},\qq
r_{k,n}^+=\max\{\l\in\ol\vk_n:\r_k(\l)=0\},\qq
\vk_n=(\l_{1,p}^{0,-},\l_{1,p}^{0,+}).
$$
Moreover, the following estimates hold true:
\begin{multline}
\lb{eeg}
E_{2,p-1}^{k,+}\le\min\{E_{2,p}^{k,-},E_{1,p-1}^{k,+}\}\le
\max\{E_{2,p}^{k,-},E_{1,p-1}^{k,+}\}
\le E_{1,p}^{k,-}\\
\le E_{1,p}^{k,+}
\le\min\{E_{1,p+1}^{k,-},E_{2,p}^{k,+}\}
\le\max\{E_{1,p+1}^{k,-},E_{2,p}^{k,+}\}\le E_{2,p+1}^{k,-},
\end{multline}
\[
\lb{mg}
E_{2,p}^{k,-}>E_{1,p-1}^{k,+}\qq \text{iff}\qq
u_k(E_{2,p}^{k,-})<0\ ;\qqq
E_{1,p+1}^{k,-}>E_{2,p}^{k,+}\qq \text{iff}\qq
u_k(E_{2,p}^{k,+})<0.
\]

\end{theorem}

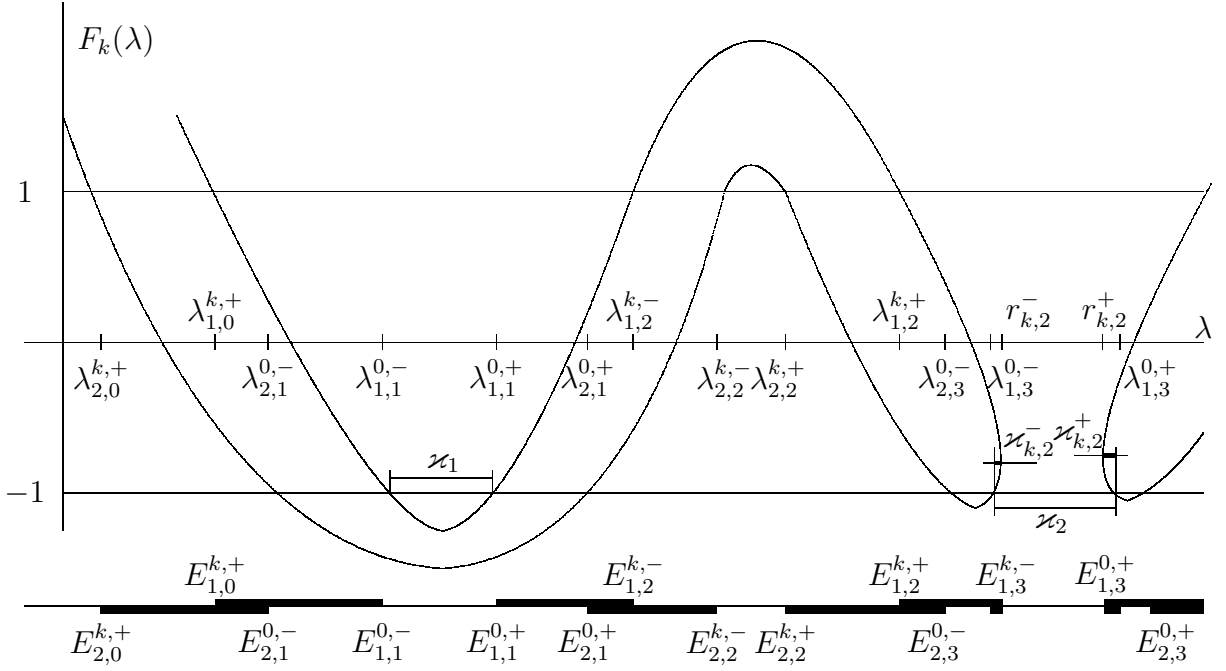
\begin{figure}
\unitlength 1.00mm \linethickness{0.2pt}
\begin{picture}(162.00,80.00)(00.00,-20.00)
\put(10.00,-5.00){\line(0,1){70.00}}
\put(10.00,40.00){\line(1,0){150.00}}
\put(05.00,20.00){\line(1,0){155.00}}
\put(10.00,00.00){\line(1,0){150.00}}
\put(160.00,22.00){\makebox(0,0)[cc]{$\l$}}
\put(05.00,00.00){\makebox(0,0)[cc]{$-1$}}
\put(05.00,40.00){\makebox(0,0)[cc]{$1$}}
\put(17.00,60.00){\makebox(0,0)[cc]{$F_k(\l)$}}
\bezier{600}(10.00,50.00)(30.00,-8.00)(60.00,-10.00)
\bezier{600}(60.00,-5.00)(50.00,-4.00)(25.00,50.00)
\bezier{400}(60.00,-5.00)(70.00,-3.00)(85.00,40.00)
\bezier{400}(85.00,40.00)(100.00,80.00)(120.00,40.00)
\bezier{400}(120.00,40.00)(140.00,02.00)(130.00,-2.00)
\bezier{400}(60.00,-10.00)(85.00,-8.00)(97.00,40.00)
\bezier{150}(97.00,40.00)(100.00,47.00)(105.00,40.00)
\bezier{400}(130.00,-2.00)(120.00,02.00)(105.00,40.00)
\bezier{400}(150.00,-1.00)(140.00,02.00)(161.00,41.00)
\bezier{250}(150.00,-1.00)(155.00,01.00)(160.00,08.00)
\put(133.50,19.00){\line(0,1){2.00}}
\put(137.00,24.00){\makebox(0,0)[cc]{$r_{k,2}^-$}}
\put(146.70,19.00){\line(0,1){2.00}}
\put(147.00,24.00){\makebox(0,0)[cc]{$r_{k,2}^+$}}
\put(15.00,19.00){\line(0,1){2.00}}
\put(15.00,15.00){\makebox(0,0)[cc]{$\l_{2,0}^{k,+}$}}
\put(30.00,19.00){\line(0,1){2.00}}
\put(30.00,24.00){\makebox(0,0)[cc]{$\l_{1,0}^{k,+}$}}
\put(85.00,19.00){\line(0,1){2.00}}
\put(85.00,24.00){\makebox(0,0)[cc]{$\l_{1,2}^{k,-}$}}
\put(120.00,19.00){\line(0,1){2.00}}
\put(120.00,24.00){\makebox(0,0)[cc]{$\l_{1,2}^{k,+}$}}
\put(96.00,19.00){\line(0,1){2.00}}
\put(97.00,15.00){\makebox(0,0)[cc]{$\l_{2,2}^{k,-}$}}
\put(105.00,19.00){\line(0,1){2.00}}
\put(104.00,15.00){\makebox(0,0)[cc]{$\l_{2,2}^{k,+}$}}
\put(37.00,19.00){\line(0,1){2.00}}
\put(37.00,15.00){\makebox(0,0)[cc]{$\l_{2,1}^{0,-}$}}
\put(52.00,19.00){\line(0,1){2.00}}
\put(52.00,15.00){\makebox(0,0)[cc]{$\l_{1,1}^{0,-}$}}
\put(67.00,19.00){\line(0,1){2.00}}
\put(67.00,15.00){\makebox(0,0)[cc]{$\l_{1,1}^{0,+}$}}
\put(79.00,19.00){\line(0,1){2.00}}
\put(79.00,15.00){\makebox(0,0)[cc]{$\l_{2,1}^{0,+}$}}
\put(126.00,19.00){\line(0,1){2.00}}
\put(126.00,15.00){\makebox(0,0)[cc]{$\l_{2,3}^{0,-}$}}
\put(132.00,19.00){\line(0,1){2.00}}
\put(135.00,15.00){\makebox(0,0)[cc]{$\l_{1,3}^{0,-}$}}
\put(149.00,19.00){\line(0,1){2.00}}
\put(153.00,15.00){\makebox(0,0)[cc]{$\l_{1,3}^{0,+}$}}
\put(53.00,00.00){\line(0,1){3.00}}
\put(66.50,00.00){\line(0,1){3.00}}
\put(60.00,04.00){\makebox(0,0)[cc]{$\vk_1$}}
\put(53.00,02.00){\line(1,0){13.50}}
\put(132.50,-03.00){\line(0,1){9.00}}
\put(148.50,-03.00){\line(0,1){9.00}}
\put(132.50,-02.00){\line(1,0){16.00}}
\put(140.00,-04.00){\makebox(0,0)[cc]{$\vk_2$}}
\put(131.00,04.00){\line(1,0){7.00}}
\put(132.50,04.20){\line(1,0){0.70}}
\put(132.50,03.80){\line(1,0){0.70}}
\put(137.00,07.00){\makebox(0,0)[cc]{$\vk_{k,2}^-$}}
\put(143.00,05.00){\line(1,0){7.00}}
\put(147.00,05.20){\line(1,0){1.50}}
\put(147.00,04.80){\line(1,0){1.50}}
\put(143.50,08.00){\makebox(0,0)[cc]{$\vk_{k,2}^+$}}
\put(05.00,-15.00){\line(1,0){155.00}}
\put(15.00,-15.90){\line(1,0){22.00}}
\put(15.00,-15.80){\line(1,0){22.00}}
\put(15.00,-15.70){\line(1,0){22.00}}
\put(15.00,-15.60){\line(1,0){22.00}}
\put(15.00,-15.50){\line(1,0){22.00}}
\put(15.00,-15.40){\line(1,0){22.00}}
\put(15.00,-15.30){\line(1,0){22.00}}
\put(15.00,-15.20){\line(1,0){22.00}}
\put(30.00,-14.80){\line(1,0){22.00}}
\put(30.00,-14.70){\line(1,0){22.00}}
\put(30.00,-14.60){\line(1,0){22.00}}
\put(30.00,-14.50){\line(1,0){22.00}}
\put(30.00,-14.40){\line(1,0){22.00}}
\put(30.00,-14.30){\line(1,0){22.00}}
\put(30.00,-14.20){\line(1,0){22.00}}
\put(30.00,-14.10){\line(1,0){22.00}}
\put(67.00,-14.10){\line(1,0){18.00}}
\put(67.00,-14.20){\line(1,0){18.00}}
\put(67.00,-14.30){\line(1,0){18.00}}
\put(67.00,-14.40){\line(1,0){18.00}}
\put(67.00,-14.50){\line(1,0){18.00}}
\put(67.00,-14.60){\line(1,0){18.00}}
\put(67.00,-14.70){\line(1,0){18.00}}
\put(67.00,-14.80){\line(1,0){18.00}}
\put(79.00,-15.90){\line(1,0){17.00}}
\put(79.00,-15.80){\line(1,0){17.00}}
\put(79.00,-15.70){\line(1,0){17.00}}
\put(79.00,-15.60){\line(1,0){17.00}}
\put(79.00,-15.50){\line(1,0){17.00}}
\put(79.00,-15.40){\line(1,0){17.00}}
\put(79.00,-15.30){\line(1,0){17.00}}
\put(79.00,-15.20){\line(1,0){17.00}}
\put(105.00,-15.90){\line(1,0){21.00}}
\put(105.00,-15.80){\line(1,0){21.00}}
\put(105.00,-15.70){\line(1,0){21.00}}
\put(105.00,-15.60){\line(1,0){21.00}}
\put(105.00,-15.50){\line(1,0){21.00}}
\put(105.00,-15.40){\line(1,0){21.00}}
\put(105.00,-15.30){\line(1,0){21.00}}
\put(105.00,-15.20){\line(1,0){21.00}}
\put(120.00,-14.10){\line(1,0){13.50}}
\put(120.00,-14.20){\line(1,0){13.50}}
\put(120.00,-14.30){\line(1,0){13.50}}
\put(120.00,-14.40){\line(1,0){13.50}}
\put(120.00,-14.50){\line(1,0){13.50}}
\put(120.00,-14.60){\line(1,0){13.50}}
\put(120.00,-14.70){\line(1,0){13.50}}
\put(120.00,-14.80){\line(1,0){13.50}}
\put(132.00,-15.90){\line(1,0){1.50}}
\put(132.00,-15.80){\line(1,0){1.50}}
\put(132.00,-15.70){\line(1,0){1.50}}
\put(132.00,-15.60){\line(1,0){1.50}}
\put(132.00,-15.50){\line(1,0){1.50}}
\put(132.00,-15.40){\line(1,0){1.50}}
\put(132.00,-15.30){\line(1,0){1.50}}
\put(132.00,-15.20){\line(1,0){1.50}}
\put(147.00,-14.10){\line(1,0){13.00}}
\put(147.00,-14.20){\line(1,0){13.00}}
\put(147.00,-14.30){\line(1,0){13.00}}
\put(147.00,-14.40){\line(1,0){13.00}}
\put(147.00,-14.50){\line(1,0){13.00}}
\put(147.00,-14.60){\line(1,0){13.00}}
\put(147.00,-14.70){\line(1,0){13.00}}
\put(147.00,-14.80){\line(1,0){13.00}}
\put(147.00,-15.90){\line(1,0){2.00}}
\put(147.00,-15.80){\line(1,0){2.00}}
\put(147.00,-15.70){\line(1,0){2.00}}
\put(147.00,-15.60){\line(1,0){2.00}}
\put(147.00,-15.50){\line(1,0){2.00}}
\put(147.00,-15.40){\line(1,0){2.00}}
\put(147.00,-15.30){\line(1,0){2.00}}
\put(147.00,-15.20){\line(1,0){2.00}}
\put(153.00,-15.90){\line(1,0){7.00}}
\put(153.00,-15.80){\line(1,0){7.00}}
\put(153.00,-15.60){\line(1,0){7.00}}
\put(153.00,-15.50){\line(1,0){7.00}}
\put(153.00,-15.40){\line(1,0){7.00}}
\put(153.00,-15.30){\line(1,0){7.00}}
\put(153.00,-15.20){\line(1,0){7.00}}
\put(153.00,-15.70){\line(1,0){7.00}}
\put(134.00,-11.00){\makebox(0,0)[cc]{$E_{1,3}^{k,-}$}}
\put(147.00,-11.00){\makebox(0,0)[cc]{$E_{1,3}^{0,+}$}}
\put(15.00,-20.00){\makebox(0,0)[cc]{$E_{2,0}^{k,+}$}}
\put(30.00,-11.00){\makebox(0,0)[cc]{$E_{1,0}^{k,+}$}}
\put(85.00,-11.00){\makebox(0,0)[cc]{$E_{1,2}^{k,-}$}}
\put(120.00,-11.00){\makebox(0,0)[cc]{$E_{1,2}^{k,+}$}}
\put(96.00,-20.00){\makebox(0,0)[cc]{$E_{2,2}^{k,-}$}}
\put(105.00,-20.00){\makebox(0,0)[cc]{$E_{2,2}^{k,+}$}}
\put(37.00,-20.00){\makebox(0,0)[cc]{$E_{2,1}^{0,-}$}}
\put(52.00,-20.00){\makebox(0,0)[cc]{$E_{1,1}^{0,-}$}}
\put(67.00,-20.00){\makebox(0,0)[cc]{$E_{1,1}^{0,+}$}}
\put(79.00,-20.00){\makebox(0,0)[cc]{$E_{2,1}^{0,+}$}}
\put(125.00,-20.00){\makebox(0,0)[cc]{$E_{2,3}^{0,-}$}}
\put(155.00,-20.00){\makebox(0,0)[cc]{$E_{2,3}^{0,+}$}}
\end{picture}
\caption{Graph of the function $F_k(\l)$ and the spectrum of $H_k$}
\lb{lfd}
\end{figure}

\no {\bf Remark.} (i) The second identity in \er{s} shows that
$E_{2,p}^{k,\pm}=E_{2,p}^{0,\pm}$
for all $(k,n)\in\Z_N\ts\N$. Here and below $p=2n-1$.

\no (ii) The last identity in \er{DeLk} gives $\r_0=9F^2$
and then $r_{0,n}^\pm=\e_n$ are zeros of $F$.

\no (iii) Let $k\ne{N\/2}$. If $v_k(\l_{1,p}^{0,\s})<0$ for some $\s=\pm$,
then $\r_k$ has at least two zeros $r_{k,n}^\pm$  in $\ol{\vk_n}$,
(see Lemma \ref{res}(iii)).
In Lemma \ref{eq} we prove that
the last identity in \er{s} for $k\ne{N\/2}$ is equivalent to
$$
E_{1,p}^{k,\pm}=\ca
\l_{1,p}^{0,\pm}\qq \text{if}\qq
F_{k,1}(r_{k,n}^\pm)=F_{k,2}(r_{k,n}^\pm)\le -1\
\text{or}\ \r_k>0\ \text{on}\ \vk_n\\
r_{k,n}^\pm\qq\ \text{if}\qq
F_{k,1}(r_{k,n}^\pm)=F_{k,2}(r_{k,n}^\pm)\in(-1,-{1\/2}]
\ac,\!\!\!\qq k\ne{N\/2}.
$$

\no (iv) Let $k=m={N\/2}\in\Z$. Then $c_m=0$ and \er{DeLk} gives
$\r_{m}=F_-^2$. Thus, $v_m=|F_-|$ and
$v_m(\l_{1,p}^{0,\pm})\ge 0$ for all $n\ge 1$, where $p=2n-1$.
The last identity in \er{s} gives
$E_{1,p}^{m,\pm}=\l_{1,p}^{0,\pm}$.

\begin{theorem} \lb{4s} 
Let $k\in\Z_N,n\ge 1,p=2n-1$.

\no (i) Let $\vk_{k,n}^-=(\l_{1,p}^{0,-},r_{k,n}^-)\ss S_{1,2n-1}^k$,
$\vk_{k,n}^+=(r_{k,n}^+,\l_{1,p}^{0,+})\ss S_{1,2n}^k$
$($i.e., $E_{1,p}^{k,\pm}=r_{k,n}^\pm)$.
Then the spectrum of $H_k$ in $\vk_{k,n}^\pm\ne\es$
has multiplicity 4.

\no (ii) If $E_{2,p}^{k,-}>E_{1,p-1}^{k,+}$
{\rm(}or $E_{1,p+1}^{k,-}>E_{2,p}^{k,+}${\rm)},
then the spectrum of $H_k$ in the interval
$(E_{1,p-1}^{k,+},E_{2,p}^{k,-})$ $=S_{1,p}^k\cap S_{2,p}^k$
{\rm(}or $(E_{2,p}^{k,+},E_{1,p+1}^{k,-})
=S_{1,p+1}^k\cap S_{2,p+1}^k${\rm)}
has multiplicity 4.

\no (iii) The spectrum $\s_{ac}(H_k)$ in all intervals, with the exception
the intervals of the statements (i), (ii), has multiplicity 2.
\end{theorem}

\no {\bf Remark} (i) Let $q\in L_{even}^2(0,1)$.
In this case $F_-=0$ (see p.8, \cite{MW}). If $k\ne{N\/2}$, then
$v_k(\l_{1,p}^{0,\pm})=-c_k^2<0$
and the last identity in \er{s} gives
$E_{1,p}^{k,\pm}=r_{k,n}^{\pm}$ for all $n\ge 1$.
If $k\ne 0$, then $u_k=-s_k^2<0$. Relations \er{mg} show that
the spectrum in each interval $S_{1,n}^k\cap S_{2,n}^k\ne\es,n\ge 1$
has multiplicity 4.

\no (ii) In Proposition \ref{pmg} we prove
that $u_k(E_{2,p}^{0,\pm})>0$ and
$v_k(\l_{1,p}^{0,\pm})>0$ for some $k, n$ and
for some specific non-even potentials.
Then relations \er{mg} give $S_{1,n}^k\cap S_{2,n}^k=\es$,
and the last identity in \er{s} yields
$E_{1,p}^{k,\pm}=\l_{1,p}^{0,\pm}$.

In order to describe gaps in the spectrum of
$H_k,H$ we need

\begin{Def}
Let $g=(\l_1,\l_2)$ be a gap in the spectrum of $H_k$ or $H$.

\no (i) If $\l_1,\l_2$ are zeros of $D_k^+$ $($or $D_k^- )$,
then $g$ is a periodic $($or antiperiodic$)$ gap.


\no (ii) If $\l_1,\l_2$ are zeros of $\r_k$,
then $g$ is a resonance gap.

\no (iii) If one of the numbers $\l_1,\l_2$ is a zero
of $D_k^-$
and other is a zero of $D_k^+$ $($or $\r_k )$,
then $g$ is a p-mix gap $($or r-mix gap$)$.


\end{Def}

\no In our armchair model there is no
a gap $(\l_1,\l_2)$, where one
of the numbers $\l_1,\l_2$
is a zero of $D_k^+$ and other is a zero of $\r_k$.

\begin{theorem}\lb{C} Let $k\in\Z_N$.
Then $\s_{ac}(H_k)=\R\sm\cup_{n\ge 0} G_{k,n}$, where the gaps
$G_{k,n}$ satisfy:
\begin{multline}
\lb{Hk0}
\wt\g_0\ss G_{k,0}=(-\iy,E_{2,0}^{k,+}),\ \
\wt\g_n\ss G_{k,4n}=(E_{2,2n}^{k,-},E_{2,2n}^{k,+}),\ \
G_{k,4n-2}=(E_{1,2n-1}^{k,-},E_{1,2n-1}^{k,+})\ss\vk_n,\\
G_{k,4n-3}=(E_{2,2n-1}^{k,-},E_{1,2n-2}^{k,+}),\qq
G_{k,4n-1}=(E_{1,2n}^{k,-},E_{2,2n-1}^{k,+}),\qq
\e_n\in[E_{1,2n-1}^{k,-},E_{1,2n-1}^{k,+}],
\end{multline}
\[
\lb{gk2}
G_{k,n}=G_{N-k,n}\ \text{all}\ \ k\in\Z_N,\qq
G_{k,4n}\ss G_{\ell,4n},\ \
G_{\ell,2n-1}\ss G_{k,2n-1}\ \
\text{all}\ \ 0\le k<\ell\le{N\/2}.
\]
Furthermore, for some $n_0\ge 1$ the gaps satisfy:

\no $G_{k,4n}$ are periodic gaps,

\no $G_{k,2n-1}$ are p-mix gaps and each
$G_{k,2n-1}=\es$  for  $k\ne 0, n\ge n_0$

\no $G_{0,4n-2}$ are  antiperiodic gaps and
$G_{0,4n-2}=\es$ for  $n\ge n_0$,

\no $G_{k,4n-2},k\not\in\{0,{N\/2}\}$ are antiperiodic,
or resonance, or r-mix gaps, and  $G_{k,4n-2}$ are resonance gaps
for $n\ge n_0$,

 \no $G_{{N\/2},4n-2},{N\/2}\in\Z$ are antiperiodic gaps.

If $q\in L^2_{even}(0,1)$,
then each $G_{k,4n-2},k\not\in\{0,{N\/2}\},n\ge 1$ is a
resonance gap.

\end{theorem}

\no{\bf Remark.} In Theorem \ref{C} and below we let the gap
$G_{k,4n-2}=\es$, if $E_{1,2n-1}^{k,-}\ge E_{1,2n-1}^{k,+}$,
and the similar relations for other gaps hold true.

Below we write  $a_n=b_n+\ell^2(n)$
for two sequences $(a_n)_1^\iy , (b_n)_1^\iy $
iff $(a_n-b_n)_1^\iy\in\ell^2$.
We describe the spectrum of $H$.

\begin{theorem} \lb{TM}
$\s_{ac}(H)=\R\sm\cup_{n\ge 0} G_{n}$, where the gaps
$G_n=\cap_{k\in\Z_N}G_{k,n}$ and $G_n$ satisfy:
\begin{multline}
\lb{TM-3}
G_0=(-\iy,E_{2,0}^+)=G_{0,0},\ \
G_{4n}=(E_{2,2n}^-,E_{2,2n}^+)=G_{0,4n},\ \
G_{4n-2}=(E_{1,2n-1}^-,E_{1,2n-1}^+)\ss\vk_n,\\
G_{4n-3}=(E_{2,2n-1}^-,E_{1,2n-2}^+)=G_{m,4n-3},\qq
G_{4n-1}=(E_{1,2n}^-,E_{2,2n-1}^+)=G_{m,4n-1},
\end{multline}
\[
\lb{TM-1}
\wt\g_{n-1}\ss G_{4n-4},
\qq\e_n\in [E_{1,2n-1}^-,E_{1,2n-1}^+]\qq\text{all}\qq n\ge 1.
\]
The gaps
$G_{4n-2}=G_{2n-1}=\es$
for all large $n\ge 1$
and the following asymptotics hold true:
\[
\lb{TM-4}
E_{2,2n}^{\pm}=E_{2,2n}^{0,\pm}=(\pi n)^2+q_0
\pm\sqrt{{2\/3}q_{sn}^2+q_{cn}^2}+{\ell^2(n)\/n}\qq \text{as}\qq
n\to \iy,
\]
where
$q_0=\int_0^1q(t)dt,
q_{sn}=\int_0^1 q(s)\sin 2\pi ns ds,\
q_{cn}=\int_0^1 q(s)\cos 2\pi ns ds.
$


\end{theorem}

There are papers about the spectral analysis of the
Schr\"odinger operator on periodic graphs and periodic nanotubes.
Molchanov and Vainberg \cite{MV} consider
 Schr\"odinger operators with $q=0$ on so-called necklace  periodic graphs.
Korotyaev and Lobanov [KL], [KL1] consider the
Schr\"odinger operator on the zigzag nanotube.
The spectrum of this operator consists of
an absolutely continuous part (intervals separated by gaps)
plus an infinite number of eigenvalues  with infinite multiplicity.
They describe all eigenfunctions with the same eigenvalue.
They define a Lyapunov function,
which is analytic on some Riemann surface. On each
sheet, the Lyapunov function has the same properties
as in the scalar case, but it has
branch points (resonances). They prove that all  resonances are
real and they determine the asymptotics of the periodic and
anti-periodic spectrum and of the resonances at high energy.
They show that there exist two types of gaps: i) stable gaps, where the endpoints are periodic and
anti-periodic eigenvalues, ii) unstable (resonance) gaps, where the
endpoints are resonances (i.e., real branch points of the Lyapunov
function).
They describe all finite gap potentials. They show that the mapping:
potential $\to$ all eigenvalues is a real analytic isomorphism
for some class of potentials.

Moreover, Korotyaev and Lobanov  [KL1] consider
magnetic Schr\"odinger operators on zigzag nanotubes.
They describe how the spectrum depends on the magnetic field.
Korotyaev [K2] considers integrated density of states and
effective masses  for zigzag nanotubes
in magnetic fields. He obtains a priori estimates of gap
lengths in terms of effective masses.
Kuchment and Post \cite{KuP} consider
the case of the zigzag, armchair  and achiral  nanotubes
with even potential $q\in L_{even}^2(0,1)$.
They show that the spectrum of the Schr\"odinger operator
(on these nanotubes), as a set, coincides with the spectrum  of the Hill
operator.

In \cite{BBKL} authors describe all eigenfunctions of $\mH$
with the same eigenvalue.
They define a Lyapunov function,
which is analytic on some Riemann surface. On each
sheet, the Lyapunov function has the same properties
as in the scalar case, but it has
branch points (resonances). They prove that there exist
non-real and real resonances.

In the present paper we describe the absolutely continuous
spectrum of $\mH$, multiplicity of the spectrum and endpoints of the spectral bands. These results are absent in \cite{KuP}. 
We show that there exist two types of gaps:
i) stable gaps, where the endpoints are periodic and
anti-periodic eigenvalues,
ii) unstable (resonance) gaps, where the
endpoints are resonances (i.e., real branch points of the Lyapunov
function).
Moreover, we determine the asymptotics of the gaps at high energy.


We present the plan of the paper.
In Sect. 2 we describe the periodic and antiperiodic eigenvalues.
In Sect. 3 we prove the main results about the spectrum of $H_k$, $\cH$.

\section{ Preliminaries }
\setcounter{equation}{0}

\begin{lemma} \lb{Rpa}
 There exists an integer $n_0>1$ such that

\no (i) The function
 $D_0^-$ given by \er{S3b}
 has exactly $4n_0$ zeros, counted with multiplicities,
 in the domain $\{\l: |\sqrt{\l}|<\pi n_0\}$
and for each $n>n_0,$ exactly two zeros,  counted with
multiplicities, in each domain
$\{\l: |\sqrt\l-\pi n-{\pi\/2}\pm\arcsin{1\/3}|<{1\/3}\}$.
There are no other zeros.

\no (ii) Each function $D_k^+,k\in\Z_N$ has
exactly $4n_0+2$ zeros,
counted with multiplicities,
 in the domain
 $\{\l: |\sqrt{\l}|<\pi n_0+{\pi\/2}\}$
and for each $n>n_0,$  exactly one simple zero
in each domain
$\{\l: |\sqrt \l-\pi n-{\pi\/2}
\pm\arcsin{\sqrt{5\pm 4c_k}\/3}|<{1\/3}\}$.
There are no other zeros.

\no (iii) Each function
$\r_k: k\not\in\{ 0, {N\/2}\}$,
has exactly $2n_0$ zeros,
counted with multiplicities,
 in the domain $\{\l: |\sqrt{\l}|<\pi n_0\}$,
and for each $n>n_0$ exactly one simple real zero
in each domain
$\{\l: |\sqrt \l-(\pi n-{\pi\/2}\pm \arcsin{s_k\/3})|
<{s_k\/3}\}$.
There are no other
zeros.

\end{lemma}

\no{\bf Proof} repeats the case of the zigzag nanotube \cite{KL}.
$\BBox$

Substituting \er{DeLk} into \er{S3b}  we
obtain for $ k\in\Z_N$
\begin{multline}
\lb{T5b1}
D_k^+=\lt((3F-1)^2-4-F_-^2\rt)\lt((3F+1)^2-4-F_-^2\rt)+16 s_k^2
=(9F^2-g_{k,1})(9F^2-g_{k,2}),\\
D_k^-=D_0^-
=\lt((3F-1)^2-F_-^2\rt)\lt((3F+1)^2-F_-^2\rt)=(9F^2-h_1)(9F^2-h_2)
\end{multline}
on $\R$,
where
\[
\lb{nec}
g_{k,\n}=5+F_-^2+(-1)^\n 2\sqrt{F_-^2+4c_k^2},\qq
h_\n=(1+(-1)^\n|F_-|)^2.
\]


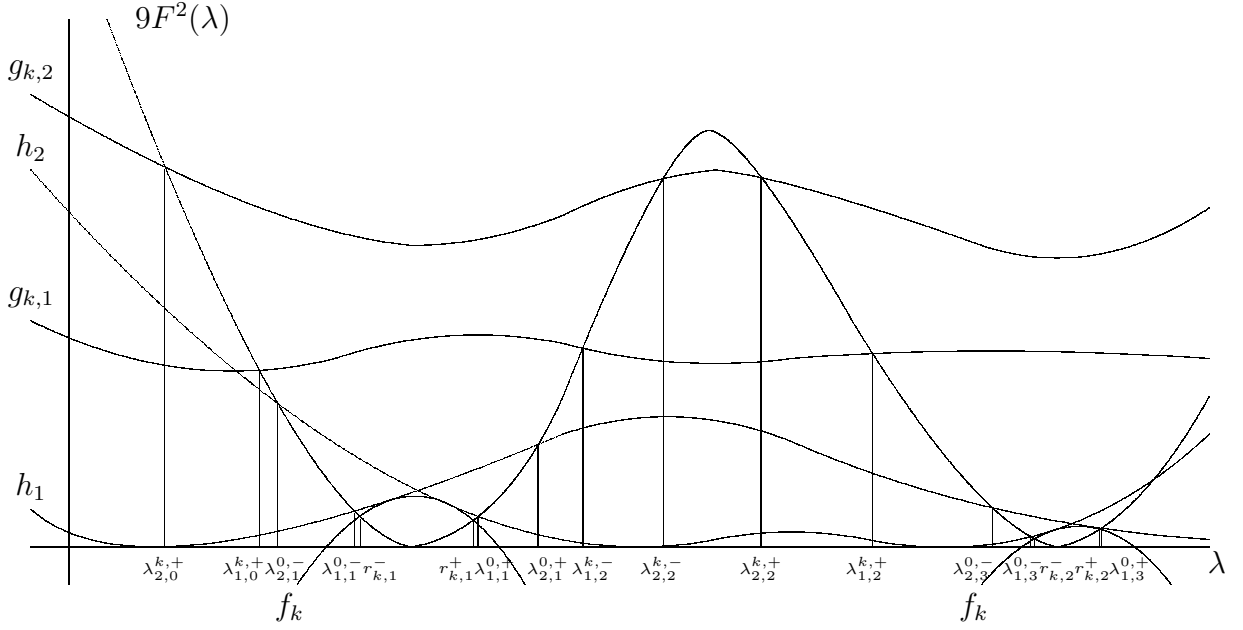
\begin{figure}
\unitlength 1.00mm \linethickness{0.1pt}
\begin{picture}(162.00,90.00)(00.00,-10.00)
\put(10.00,-5.00){\line(0,1){75.00}}
\put(05.00,00.00){\line(1,0){155.00}}
\put(161.00,-2.00){\makebox(0,0)[cc]{$\l$}}
\put(25.00,70.00){\makebox(0,0)[cc]{$9F^2(\l)$}}
\bezier{600}(15.00,70.00)(40.00,3.00)(55.00,0.00)
\bezier{600}(55.00,0.00)(67.00,2.00)(75.00,20.00)
\bezier{600}(75.00,20.00)(90.00,58.00)(95.00,55.00)
\bezier{600}(95.00,55.00)(100.00,53.00)(113.00,30.00)
\bezier{600}(113.00,30.00)(130.00,02.00)(140.00,00.00)
\bezier{600}(140.00,00.00)(150.00,02.00)(160.00,20.00)
\bezier{600}(5.00,60.00)(35.00,42.00)(55.00,40.00)
\bezier{600}(55.00,40.00)(65.00,40.00)(75.00,44.00)
\bezier{600}(75.00,44.00)(85.00,49.00)(95.00,50.00)
\bezier{600}(95.00,50.00)(105.00,49.00)(130.00,40.00)
\bezier{600}(130.00,40.00)(145.00,35.00)(160.00,45.00)
\put(05.00,63.00){\makebox(0,0)[cc]{$g_{k,2}$}}
\bezier{600}(5.00,30.00)(25.00,20.00)(45.00,25.00)
\bezier{600}(45.00,25.00)(60.00,30.00)(75.00,27.00)
\bezier{600}(75.00,27.00)(90.00,23.00)(105.00,25.00)
\bezier{600}(105.00,25.00)(130.00,27.00)(160.00,25.00)
\put(05.00,33.00){\makebox(0,0)[cc]{$g_{k,1}$}}
\bezier{600}(5.00,50.00)(55.00,-7.00)(95.00,1.00)
\bezier{600}(95.00,1.00)(105.00,3.00)(115.00,1.00)
\bezier{600}(115.00,1.00)(140.00,-4.00)(160.00,15.00)
\put(05.00,53.00){\makebox(0,0)[cc]{$h_{2}$}}
\bezier{600}(5.00,5.00)(20.00,-8.60)(75.00,15.00)
\bezier{600}(75.00,15.00)(90.00,20.00)(105.00,14.00)
\bezier{600}(105.00,14.00)(130.00,3.00)(160.00,1.00)
\put(05.00,8.00){\makebox(0,0)[cc]{$h_{1}$}}
\bezier{600}(40.00,-5.00)(56.00,18.50)(70.00,-5.00)
\bezier{600}(130.00,-5.00)(143.20,10.50)(155.00,-5.00)
\put(39.00,-8.00){\makebox(0,0)[cc]{$f_{k}$}}
\put(129.00,-8.00){\makebox(0,0)[cc]{$f_{k}$}}
\put(22.52,0.00){\line(0,1){50.00}}
\put(22.50,-3.00){\makebox(0,0)[cc]{\tiny $\l_{2,0}^{k,+}$}}
\put(35.00,0.00){\line(0,1){23.50}}
\put(33.00,-3.00){\makebox(0,0)[cc]{\tiny $\l_{1,0}^{k,+}$}}
\put(37.40,0.00){\line(0,1){19.00}}
\put(38.50,-3.00){\makebox(0,0)[cc]{\tiny $\l_{2,1}^{0,-}$}}
\put(47.50,0.00){\line(0,1){4.50}}
\put(46.00,-3.00){\makebox(0,0)[cc]{\tiny $\l_{1,1}^{0,-}$}}
\put(48.30,0.00){\line(0,1){3.70}}
\put(51.00,-3.00){\makebox(0,0)[cc]{\tiny $r_{k,1}^{-}$}}
\put(63.20,0.00){\line(0,1){3.30}}
\put(61.00,-3.00){\makebox(0,0)[cc]{\tiny $r_{k,1}^{+}$}}
\put(63.80,0.00){\line(0,1){4.10}}
\put(66.00,-3.00){\makebox(0,0)[cc]{\tiny $\l_{1,1}^{0,+}$}}
\put(71.70,0.00){\line(0,1){13.50}}
\put(73.00,-3.00){\makebox(0,0)[cc]{\tiny $\l_{2,1}^{0,+}$}}
\put(77.60,0.00){\line(0,1){26.30}}
\put(79.00,-3.00){\makebox(0,0)[cc]{\tiny $\l_{1,2}^{k,-}$}}
\put(88.15,0.00){\line(0,1){48.80}}
\put(88.00,-3.00){\makebox(0,0)[cc]{\tiny $\l_{2,2}^{k,-}$}}
\put(101.00,0.00){\line(0,1){48.90}}
\put(101.00,-3.00){\makebox(0,0)[cc]{\tiny $\l_{2,2}^{k,+}$}}
\put(115.70,0.00){\line(0,1){25.50}}
\put(115.00,-3.00){\makebox(0,0)[cc]{\tiny $\l_{1,2}^{k,+}$}}
\put(131.50,0.00){\line(0,1){5.00}}
\put(129.00,-3.00){\makebox(0,0)[cc]{\tiny $\l_{2,3}^{0,-}$}}
\put(136.45,0.00){\line(0,1){1.50}}
\put(135.50,-3.00){\makebox(0,0)[cc]{\tiny $\l_{1,3}^{0,-}$}}
\put(137.00,0.00){\line(0,1){1.00}}
\put(140.00,-3.00){\makebox(0,0)[cc]{\tiny $r_{k,2}^{-}$}}
\put(145.50,0.00){\line(0,1){2.30}}
\put(144.60,-3.00){\makebox(0,0)[cc]{\tiny $r_{k,2}^{+}$}}
\put(145.80,0.00){\line(0,1){2.50}}
\put(149.50,-3.00){\makebox(0,0)[cc]{\tiny $\l_{1,3}^{0,+}$}}
\end{picture}
\caption{Functions $9F^2,g_k^\pm,h^\pm$ and $f_k$}
\lb{lf}
\end{figure}

\begin{lemma}\label{T4n}
(i) For all  $(\n,k,n)\in\N_2\ts\Z_N\ts\N$
the periodic and antiperiodic eigenvalues
satisfy
\[
\lb{sym}
\l_{\n,n-1}^{k,\pm}=\l_{\n,n-1}^{N-k,\pm},\qq
\l_{\n,2n-1}^{k,\pm}=\l_{\n,2n-1}^{0,\pm},\qq
\g_{\n,2n-1}^k=\g_{\n,2n-1}^0,\qq \g_{1,2n-1}^0=\vk_n,
\]
\begin{multline}
\lb{b1}
\wt\l_{n-1}^+\le\l_{2,p-1}^{k,+}
\le\min\{\l_{2,p}^{0,-},\l_{1,p-1}^{k,+}\}
\le\max\{\l_{2,p}^{0,-},\l_{1,p-1}^{k,+}\}
\le\l_{1,p}^{0,-}\\ \le\e_n\le\l_{1,p}^{0,+}
\le\min\{\l_{1,p+1}^{k,-},\l_{2,p}^{0,+}\}
\le\max\{\l_{1,p+1}^{k,-},\l_{2,p}^{0,+}\}
\le\l_{2,p+1}^{k,-}\le\wt\l_{n}^-,\qq p=2n-1,
\end{multline}
\[
\lb{cg}
\l_{2,p}^{0,-}<\l_{1,p-1}^{k,+}\ \Leftrightarrow\
u_k(\l_{2,p}^{0,-})>0\ ;
\qqq\l_{2,p}^{0,+}>\l_{1,p+1}^{k,-}\ \Leftrightarrow\
u_k(\l_{2,p}^{0,+})>0,
\]
\[
\lb{s1a}
\cup_{n\ge 1}
\g_{\n,2n-1}^0
=\{\l\in\R:9F^2(\l)<h_\n(\l)\},\qq
\cup_{n\ge 0}\g_{\n,2n}^k
=\{\l\in\R:9F^2(\l)>g_{k,\n}(\l)\},
\]
where
$$
\g_{\n,0}^k=(-\iy,\l_{\n,0}^{k,+}),\qq
\g_{\n,n}^k=(\l_{\n,n}^{k,-},\l_{\n,n}^{k,+}),\qq
(\n,n,k)\in\N_2\ts\N\ts\Z_N.
$$

\no (ii) If $0\le k<\ell\le{N\/2}$ and $n\ge 1$, then
\[
\lb{F2}
\l_{2,2n-2}^{k,+}<\l_{2,2n-2}^{\ell,+},\qq
\l_{1,2n-2}^{\ell,+}<\l_{1,2n-2}^{k,+},\qq
\l_{1,2n}^{k,-}<\l_{1,2n}^{\ell,-},\qq
\l_{2,2n}^{\ell,-}<\l_{2,2n}^{k,-}.
\]

\end{lemma}

\no {\bf Proof.} (i) The periodic eigenvalues are zeros $D_k^+$.
Using \er{T5b1} and the identities $s_{N-k}=s_k$ we obtain
the first identity in \er{sym} for the periodic eigenvalues.
The antiperiodic eigenvalues are zeros $D_k^-$.
Using \er{T5b1} and the definitions of $\g_{\n,n}^k,\vk_n$ we obtain
the other identities in \er{sym}.

Identities \er{nec} give
$
g_{k,1}-h_1=2(|F_-|+2-\sqrt{F_-^2+4c_k^2})\ge 0
$
on $\R$. Then we obtain
\[
\lb{efs}
h_1
\le\min\{h_2,g_{k,1}\}\le\max\{h_2,g_{k,1}\}
\le g_{k,2}.
\]
Identities \er{Ho} give $F^2-F_-^2=\vt_1\vp_1'=1+\vt_1'\vp_1$.
Then $F_-^2(\m_n)=F^2(\m_n)-1$.
The last identity and $F^2(\m_n)\ge 1$ imply
$$
9F^2(\m_n)-g_{k,2}(\m_n)
=2\lt(\lt(2\sqrt{F^2(\m_n)-1+4c_k^2}-{1\/4}\rt)^2
+{31\/16}-16c_k^2\rt)\ge 0,
$$
which yields $g_{k,2}(\m_n)\le 9F^2(\m_n)$.
Estimates \er{efs} and
the properties of the function $F$ provide
that each of the functions
$9F^2-g_{k,\n}$, $9F^2-h_\n,\n=1,2$ has at least
one zero in each of the
intervals $(-\iy,\e_1]$, $[\e_n,\m_n]$, $[\m_n,\e_{n+1}]$, $n\ge 1$.
Moreover, Lemma \ref{Rpa} shows that each of these functions
has exactly one zero in each of these intervals.
Then the properties of the function $F$
and estimates \er{efs}
yield \er{b1}.

Identities \er{nec} give
$
g_{k,1}-h_2=2(s_k^2-|F_-|+\sqrt{s_k^4+4c_k^2}
-\sqrt{F_-^2+4c_k^2}).
$
For fixed $\l$ we obtain $g_{k,1}(\l)>h_2(\l)$ iff $|F_-(\l)|<s_k^2$.
Let $\l=\l_{2,2n-1}^{0,-}$. Estimates \er{b1} show that
$\wt\l_{n-1}^+\le\min\{\l_{2,2n-2}^{k,-},\l_{2,2n-1}^{0,-}\}
\le\max\{\l_{2,2n-2}^{k,-},\l_{2,2n-1}^{0,-}\}\le\e_n$.
Since $(F^2)'<0$ on $(\wt\l_{n-1}^+,\e_n)$, we deduce that
$\l_{2,2n-1}^{0,-}>\l_{2,2n-2}^{k,-}$ iff $|F_-(\l_{2,2n-1}^{0,-})|<s_k^2$,
which yields the first equivalence in \er{cg}.
The proof of the second equivalence is similar.

 Estimates \er{b1}, \er{efs} and the properties
of $F$ yield \er{s1a}.

\no (ii) Identities \er{nec} give
$g_{k,2}<g_{\ell,2},\ g_{\ell,1}<g_{k,1}$.
The properties of $F$ yield \er{F2}.
$\BBox$

\section{ Proof of Theorems \ref{Tk}--\ref{TM}}
\setcounter{equation}{0}

Let $R_k=\{\l\in\R:\r_{k}(\l)>0\}, k\ne{N\/2}$.
In Lemmas \ref{sp}-\ref{spe} we describe the set
$\s_{k,\n}=\{\l\in\R:F_{k,\n}(\l)\in[-1,1]\}$ in terms of  $F$.

\begin{lemma} \lb{sp}
 For all $k\in\Z_N$ and $\l\in R_k$
the following identities hold true:
\[
\lb{F1+}
F_{k,\n}(\l)<1\qq \text{iff}\qq
9F^2(\l)<g_{k,\n}(\l),\qq \n=1,2,
\]
\[
\lb{g1}
F_{k,1}(\l)>-1\qq \text{iff}\qq
\{9F^2(\l)>h_1(\l)\ \text{or}\
|F_-(\l)|<c_k^2\},
\]
\[
\lb{g2}
F_{k,2}(\l)>-1\qq \text{iff}\qq
\lt\{9F^2(\l)>h_2(\l)\ \text{or}\
\{9F^2(\l) < h_1(\l)\ \text{and}\
|F_-(\l)|<c_k^2\}\rt\}.
\]
\end{lemma}

\no {\bf Proof.}
If $k=m={N\/2}\in\Z$, identities \er{DeLk}, \er{nec} give
$
F_{m,\n}-1={1\/2}(9F^2-g_{m,\n}),\
F_{m,\n}+1={1\/2}(9F^2-h_\n),
$
which yields \er{F1+}-\er{g2} for $k=m$.

Let $k\ne{N\/2}$. We rewrite the functions $\r_k,F_{k,\n}-1$
in the form
$$
\r_k=(9F^2-g_{k,\n})c_k^2+(\sqrt{F_-^2+4c_k^2}+(-1)^\n c_k^2)^2,
$$
$$
F_{k,\n}-1={1\/2c_k^2}\lt(\sqrt{\r_k}-(-1)^\n c_k^2
-\sqrt{F_-^2+4c_k^2}\rt)\lt(\sqrt{\r_k}-(-1)^\n c_k^2
+\sqrt{F_-^2+4c_k^2}\rt).
$$
These identities yields \er{F1+}. 
We rewrite the functions $\r_k,F_{k,1}+1$ in the form
\[
\lb{r4}
\r_k=(9F^2-h_1)c_k^2+(|F_-|-c_k^2)^2,\ \
F_{k,1}+1={1\/2c_k^2}\lt(\sqrt{\r_k}+c_k^2
-|F_-|\rt)\lt(\sqrt{\r_k}+c_k^2
+|F_-|\rt).
\]
These identities imply \er{g1}. 
We rewrite the functions $\r_k,F_{k,2}+1$ in the form
$$
\r_k=(9F^2-h_2)c_k^2+(|F_-|+c_k^2)^2,\qq
F_{k,2}+1={1\/2c_k^2}\lt(\sqrt{\r_k}-c_k^2
-|F_-|\rt)\lt(\sqrt{\r_k}-c_k^2+|F_-|\rt).
$$
These identities and the first identity in \er{r4} give \er{g2}.
\BBox

Now we describe the zeros of $\r_k$ and the functions $F_{k,\n}$
on the interval $\vk_n$.
Recall the intervals
$\vk_{n,k}^-=(\l_{1,2n-1}^{0,-},r_{k,n}^-)$,
$\vk_{n,k}^+=(r_{k,n}^+,\l_{1,2n-1}^{0,+}),n\ge 1$ (see Theorem \ref{4s}).

\begin{lemma} \lb{res}
 Let $k\ne{N\/2}$.
Then

\no (i) The following relation holds true:
\[
\lb{r>0}
\R\sm\ol{R_k}\ss\cup_{n\ge 1}\vk_n,\qq \text{where}\qq
\vk_n=(\l_{1,2n-1}^{0,-},\l_{1,2n-1}^{0,+}).
\]
All real zeros of $\r_k$ belong to the set $\cup_{n\ge 1}\ol\vk_n$.
Each interval $\ol\vk_n,n\ge 1$ contains
even number $\ge 0$ of zeros of $\r_k$, counted with multiplicities.

\no (ii) Let $\vk_n\not\ss R_k$ for some $n\ge 1$.
Then for each $\s=\pm$ the following
relations hold true:
\[
\lb{sFc}
\sign v_k=\const\ \ \text{on\ each}\ \
\vk_{k,n}^\s\ss R_k,
\]
\[
\lb{rg1}
\text{if}\qq v_k(\l)>0\ \text{for\ some}\ \l\in \vk_{k,n}^\s,\qq
\text{then}\qq F_{k,2}<F_{k,1}<-1\ \text{on}\   \vk_{k,n}^\s,
\]
\[
\lb{rg2}
\text{if}\qq v_k(\l)<0\ \text{for\ some}\ \l\in \vk_{k,n}^\s,\qq
\text{then}\qq -1<F_{k,2}<F_{k,1}\ \text{on}\   \vk_{k,n}^\s.
\]
If $\vk_n\ss R_k$ for some $n\ge 1$, then
\[
\lb{un}
v_k>0\qq \text{and}\qq F_{k,2}<F_{k,1}<-1\qq\text{on}\qq  \vk_n.
\]

\no (iii) If $v_k(\l)<0$ for some $\l\in\ol{\vk_n},n\ge 1$, then $\r_k$
has even number $\ge 2$ of zeros on $\ol{\vk_n}$.
\end{lemma}

\no {\bf Proof.}
(i)
The last identity in \er{DeLk} gives
\[
\lb{f}
\r_k=c_k^2(9F^2-f_k),\qq k\ne{N\/2},\qqq \text{where}\qq
f_k=s_k^2\lt(1-{F_-^2\/c_k^2}\rt).
\]
Identities \er{f} show that zeros of $\r_k$ (resonances)
are zeros of $9F^2-f_k$. Identities \er{nec}, \er{f} give
$h_1-f_k=(c_k-{F_-^2\/c_k})^2\ge 0$.
Identities \er{b1} and the properties of the function $F$
imply that real zeros of $\r_k$
belong to the set $\cup_{n\ge 1}\ol\vk_n$ and
\er{r>0} holds.
The function $\r_k$ has even number of zeros in
$\ol\vk_n$, since $\r_k\ge 0$
at the points $\l_{1,2n-1}^{0,-},\l_{1,2n-1}^{0,+}$.

\no (ii) Identity \er{s1a} show that $9F^2<h_1$ on $\vk_n$.
Using the first identity in \er{r4} we conclude that if
$|F_-(\l)|=c_k^2$ for $\l\in \vk_n$, then $\r_k(\l)<0$.
We obtain \er{sFc}, since $\r_k>0$ on $\vk_{k,n}^\pm$.

If $|F_-(\l)|>c_k^2$ for some $\l\in \vk_{k,n}^\pm$, then
\er{sFc} show $|F_-(\l)|>c_k^2$ for all $\l\in \vk_{k,n}^\pm$.
Recall that $9F^2<h_1$ on $\vk_n^\pm$.
Then \er{g1} gives $F_{k,1}<-1$ on $\vk_{k,n}^\pm$,
which yields \er{rg1}.

If $|F_-(\l)|<c_k^2$ for some $\l\in \vk_{k,n}^\pm$, then
\er{sFc} provide $|F_-(\l)|<c_k^2$ for all $\l\in \vk_{k,n}^\pm$.
Using $9F^2<h_1$ on $\vk_n^\pm$ again \er{g2} gives
$-1<F_{k,2}$ on $\vk_{k,n}^\pm$. We obtain \er{rg2}.

Suppose that $\vk_n\in R_k$, i.e. $\r_k>0$ on $\vk_n$.
Estimates \er{b1} yield
$\e_n\in \vk_n$, hence $\r_k(\e_n)>0$.
Note that $\r_0(\e_n)=0$, hence the condition $\vk_n\in R_0$
is not fulfilled for all $n$.
In this reason we assume below $k\ne 0$.

The last identity in
\er{DeLk} gives $F_-^2(\e_n)=s_k^{-2}\r_k(\e_n)+c_k^2>c_k^2$, which yield
$|F_-(\e_n)|>c_k^2$.
Using $9F^2<h_1$ on $\vk_n$ and the first identity in \er{r4}
again we conclude that if
$|F_-(\l)|=c_k^2$ for $\l\in \vk_n$, then $\r_k(\l)<0$, which yields
$\sign(|F_-|-c_k^2)=\const$ on $\vk_n$.
Thus $|F_-(\l)|>c_k^2$
for all $\l\in \vk_n$ and we have the first estimate in \er{un}.
Relation \er{g1} gives
$F_{k,1}<-1$ on $\vk_n$, which yields the second estimate in
\er{un}.

\no (iii) Using relation \er{un} we deduce that
if $v_k(\l)<0$ for some
$\l\in\ol{\vk_n}=[\l_{1,2n-1}^{0,-},\l_{1,2n-1}^{0,+}]$,
then $\vk_n\not\ss R_k$.
Hence there exists $\l\in\vk_n$ such that $\r_k(\l)<0$.
On the other hand $\r_k(\l_{1,2n-1}^{0,\pm})\ge 0$,
which yields the needed statement.
$\BBox$

\no For each $(k,\n)\in\Z_N\ts\N_2$ we introduce the sets
\[
\lb{ss}
\gS_{k,\n}=\bigcup_{n\ge 1}\lt([\l_{\n,2n-2}^{k,+},\l_{\n,2n-1}^{0,-}]
\cup[\l_{\n,2n-1}^{0,+},\l_{\n,2n}^{k,-}]
\rt),\qq
\gS^R_{k}=\bigcup_{\s=\pm,n\in N_k^\s}\!\!\!\!
\ol{\vk_{k,n}^\s},
\]
where $N_k^\pm=\{n\in\N:v_k(\l_{1,2n-1}^{0,\pm})<0\}.$
The set $\gS_{k,\n}$ is a part of $\s_{k,\n}$, where the periodic
and antiperiodic eigenvalues are endpoints of bands.
The set $\gS^R_{k}$ is an "unstable" part of $\s_{k,\n}$.

\begin{lemma} \lb{spe}
For each $(\n,k)\in\N_2\ts\Z_N$
the following identities hold true:
\[
\lb{cgS}
\gS_{k,\n}=\{\l\in\R:h_\n(\l)\le 9F^2(\l)\le g_{k,\n}(\l)\},
\]
\[
\lb{g2b}
\gS^R_{k}=\{\l\in R_k:9F^2(\l)\le h_1(\l)\ \text{and}\
v_k(\l)\le 0\},\qq k\ne{N\/2},\qq \text{and}\qq
\gS^R_{N\/2}=\es,
\]
\[
\lb{F21}
\s_{k,\n}=\gS_{k,\n}\cup\gS^R_{k}.
\]

\end{lemma}

\no {\bf Proof.} Identities \er{s1a} give \er{cgS}.
If $m={N\/2}\in\Z$, then $c_m=0$,
$v_m(\l)=|F_-(\l)|\ge 0$ and
$N_k^\pm=\es$, which yields $\gS^R_{m}=\es$.
The first identity in \er{s1a} gives
$\cup_{n\ge 1}\vk_{n}=\{\l\in\R:9F^2(\l)<h_1(\l)\}$.
Then \er{rg1}-\er{un}  provide
\er{g2b} for $k\ne{N\/2}$.
Identities \er{F1+}-\er{g2} yield \er{F21}.
$\BBox$

We prove our main results.

\no{\bf Proof of Theorem \ref{Tk}.}
Identities \er{ss}-\er{F21} give
$$
\s_{k,1}
=\cup_{n\ge 1}(S_{1,2n-1}^k
\cup S_{1,2n}^k),
\qq
\s_{k,2}
=(\cup_{\s=\pm,n\in N_k^\s}
\ol{\vk_{k,n}^\s})\cup
(\cup_{n\ge 1}(S_{2,2n-1}^k
\cup S_{2,2n}^k
)),
$$
where
$S_{1,2n-1}^k=[\l_{1,2n-2}^{k,+},\l_{1,2n-1}^{0,-}]\cup\ol{\vk_{k,n}^-}$
and
$S_{1,2n}^k=\ol{\vk_{k,n}^+}\cup[\l_{\n,2n-1}^{0,+},\l_{\n,2n}^{k,-}]$.
Then \er{s} holds true.
Using \er{T3-1} and
$\ol{\vk_{k,n}^-}\ss S_{1,2n-1}^k,\ol{\vk_{k,n}^+}\ss S_{1,2n}^k$,
we obtain
$\s_{ac}(H_k)=\s_{k,1}
\cup\s_{k,2}=\cup_{\n\in\N_2,n\ge 1}S_{n,\n}^k.$
Estimates \er{b1} give \er{eeg}.
Relations \er{cg} provide \er{mg}.$\BBox$

\no{\bf Proof of Theorem \ref{4s}.}
The last identity in \er{s} together with \er{rg2} imply
$E_{1,p}^{k,\pm}=r_{k,n}^\pm$ iff
$-1<F_{k,2}<F_{k,1}$ on $\vk_{k,n}^\pm$. Identity \er{ss} shows
that $\vk_{k,n}^\pm\ss\gS_k^R$. Identity \er{g2b} yields $9F^2\le h_1$
on $\vk_{k,n}^\pm$. Relations \er{g1} give
$F_{k,2}<F_{k,1}<1$ on $\vk_{k,n}^\pm$.
Then the spectrum in $\vk_{k,n}^\pm$ has multiplicity 4.
Suppose that $E_{1,p}^{k,\pm}\ne r_{k,n}^\pm$.
Then the last identity in \er{s} show $v_k(\l_{1,2n-1}^{0,\pm})\ge 0$.
Relation \er{rg1} yields $F_{k,2}<F_{k,1}<-1$ on $\vk_{k,n}^\pm$.
Hence the interval $\vk_{k,n}^\pm$ lies in a gap of $H_k$.

Using \er{s} we rewrite $\gS_{k,\n}$ (see \er{ss}) in the form
$
\gS_{k,\n}=\bigcup_{n\ge 1}([E_{\n,2n-2}^{k,+},E_{\n,2n-1}^{0,-}]
\cup[E_{\n,2n-1}^{0,+},E_{\n,2n}^{k,-}]
).
$
Estimates \er{eeg} show that
$$
\wt\gS_k=\gS_{k,1}\cap\gS_{k,2}=
\lt(\cup_{n:E_{2,p}^{k,-}<E_{1,p-1}^{k,+}}
[E_{2,p}^{k,-},E_{1,p-1}^{k,+}]\rt)\bigcup
\lt(\cup_{n:E_{2,p}^{k,+}<E_{1,p+1}^{k,-}}
[E_{2,p}^{k,+},E_{1,p+1}^{k,-}]\rt).
$$
Identity \er{F21} give
$F_{k,\n}\in[-1,1]$ on $\gS_{k,\n}$.
Then $F_{k,\n}\in[-1,1]$ for $\n=1,2$ on $\wt\gS_k$.
Hence the spectrum on this
set has multiplicity 4.
The spectrum on $(\gS_{k,1}\cup\gS_{k,2})\sm\wt\gS_k$ has multiplicity 2.
%
$\BBox$

We need the following well known asymptotics (see, for example, \cite{K})
\begin{multline}
\lb{as}
F(\l)= \cos\sqrt{\l}+ {q_0\sin
\sqrt{\l}\/2\sqrt{\l}}+{O(e^{|\Im\sqrt{\l}|})\/|\l|},\\
F_-(\l)=-{1\/2\sqrt \l}\int_0^1\sin \sqrt \l(1-2t)q(t)dt
+{O(e^{|\Im\sqrt{\l}|})\/|\l|},\qq |\l|\to \iy.
\end{multline}

\no{\bf Proof of Theorem \ref{C}.}
Estimates \er{b1} show that the
intervals
$G_{\n,n}^k=(E_{\n,n}^{k,-},E_{\n,n}^{k,+})$, $(\n,n)\in\N_2\ts\N$,
satisfy:
$$
G_{1,0}^k\cap G_{2,m}^k=\es\  \text{for}\ m\not\in\{0,1\},\qq
 G_{1,2n-1}^k\cap G_{2,m}^k=\es\  \text{for}\ m\ne 2n-1,
$$
$$
 G_{1,2n}^k\cap G_{2,m}^k=\es\  \text{for}\ m\not\in\{2n-1,2n,2n+1\}.
$$
Then the gaps $G_{k,n},n\ge 0$ in the spectrum $H_k$
are given by
$$
G_{k,0}= G_{1,0}^k\cap G_{2,0}^k,\
G_{k,2n}= G_{1,n}^k\cap G_{2,n}^k,\
G_{k,4n-3}= G_{1,2n-2}^k\cap G_{2,2n-1}^k,\
G_{k,4n-1}= G_{1,2n}^k\cap G_{2,2n-1}^k,
$$
$ n\ge 1$, which yields all identities in \er{Hk0}.
Estimates \er{b1} give all inclusions in \er{Hk0}.
Lemma \ref{T4n} and relations \er{r>0}, \er{F2}
give \er{gk2}.
Identities \er{s} show that
$G_{k,4n}$ are periodic gaps,
$G_{k,2n-1}$ are p-mix gaps and
$G_{k,4n-2}$ are antiperiodic, or resonance,
or r-mix gaps.
Asymptotics \er{as}
and estimates \er{mg} give that
$E_{k,4n-3}^->E_{k,4n-3}^+$ and
$E_{k,4n-1}^->E_{k,4n-1}^+$ for $k\ne 0$ and large $n>1$.
Hence $G_{k,2n-1}=\es$ for such $k,n$.

Recall that $r_{0,n}^-=r_{0,n}^+=\e_n$.
Identities \er{s} give
$
E_{1,p}^{0,\pm}=\ca
\l_{1,p}^{0,\pm}\qq \text{if}\qq
v_0(\l_{1,p}^{0,\pm})\ge 0\\
\e_n\qq\ \text{if}\qq
v_0(\l_{1,p}^{0,\pm})<0
\ac\!\!\!
$.
Hence
$G_{0,4n-2}$ are antiperiodic gaps or $G_{0,4n-2}=\es$.
Moreover, $E_{1,p}^{0,-}=E_{1,p}^{0,+}=\e_n$
for all large $n\ge 1$. Hence $G_{0,4n-2}=\es$ for large $n\ge 1$.

Since $c_m=0,m={N\/2}\in\Z$, identities \er{s} provide
$E_{1,p}^{m,\pm}=\l_{1,p}^{m,\pm}$.
Hence $G_{m,4n-2}$ are antiperiodic gaps.

For $k\ne \{0,{N\/2}\}$ identities \er{s} give
$G_{k,4n-2}$ are antiperiodic, or resonance, or r-mix gaps,
and asymptotics \er{as} show that $G_{k,4n-2}$ are
resonance gaps.

If $q\in L_{even}^2(0,1)$, then $F_-=0$ and $v_k<0,k\ne{N\/2}$.
Identities \er{s} show that
$E_{1,p}^{k,\pm}=r_{k,n}^\pm$ in this case.
The last identity in \er{DeLk} yield $\r_k=(9F^2-s_k^2)c_k^2$.
Properties on the function $F$ show that $r_{k,n}^-<r_{k,n}^+$
for $k\not\in{0,{N\/2}}$ and all $n\ge 1$.
Then  $G_{k,4n-2}=(r_{k,n}^-,r_{k,n}^+),n\ge 1$ are resonance gaps.$\BBox$


\no{\bf Proof of Theorem \ref{TM}.}
Recall that the operator $\mH$ is unitarily equivalent to
$H=\os_1^N H_k$.
Relations \er{gk2} provide
$\s_{ac}(H)=\R\sm\cup_{n\ge 0} G_{n}$, where gap
$G_n=\cap_{k\in\Z_N}G_{k,n}$.
The second relations in \er{gk2} show
$G_{4n}
=G_{0,4n},n\ge 0$.
The third relations in \er{gk2} imply $G_{2n-1}=G_{m,2n-1},n\ge 1$.
the relations $G_{k,4n-2}\ss\vk_n$ give $G_{4n-2}\ss\vk_n$.
Thus, we have proved all relations in \er{TM-3}.
The relations $\wt\g_{n-1}\ss G_{k,4n},
\e_n\in [E_{1,2n-1}^{k,-},E_{1,2n-1}^{k,+}]$
give the corresponding relations
$\wt\g_{n}\ss G_{4n},
\e_n\in [E_{1,2n-1}^-,E_{1,2n-1}^+]$,
which yields \er{TM-1}.
By Theorem \ref{C}, $G_{m,2n-1}=G_{0,4n-2}=\es$
for large $n>1$, which implies $G_{4n-2}=G_{2n-1}=\es$
for large $n\ge 1$.

In order to prove asymptotics \er{TM-4} we assume that $\int_0^1q(t)dt=0$.
Identities \er{T5b1} and Lemma \ref{Rpa} (ii)
show that $\l_{2,2n}^{0,\pm}$ are zeros of the equation
\[
\lb{L54c}
3F(\l)=(-1)^n(1+\sqrt{F_-^2(\l)+4}),
\]
and $\l_{2,2n}^{0,\pm}=(\pi n+\ve_{n}^\pm)^2$,
where
 $|\ve_{n}^\pm|\le {1\/3}$ for large $n$.
Let $\l=\l_{2,2n}^{0,\pm}$ and $\ve=\ve_{n}^\pm$.
Asymptotics \er{as}
give
$$
F(\l)=(-1)^n+O(\ve^2)+O(n^{-2}),\qq F_-(\l)=O(n^{-1}).
$$
Substituting these asymptotics into \er{L54c} we get $\ve=O(n^{-1})$.
Using the standard calculations (see \cite{K}) we obtain
$$
F(\l)=(-1)^n\lt(1+{q_{sn}^2+q_{cn}^2\/2(2\pi
n)^2}-{\ve^2\/2}\rt)+{\ell^2(n)\/n^3},\qq
F_-(\l)={(-1)^nq_{sn}\/2\pi n}+{\ell^2(n)\/n^2}.
$$
Substituting the last asymptotics into \er{L54c} we obtain
$
\ve^2={{2\/3}q_{sn}^2+q_{cn}^2\/2(2\pi n)^2}+{\ell^2(n)\/n^3},
$
which yields
$
\l_{2,2n}^{0,\pm}=(\pi n)^2\pm\sqrt{{2\/3}q_{sn}^2+q_{cn}^2}+{\ell^2(n)\/n}.
$
Identity \er{TM-3} yield
$\l_{2,2n}^{0,\pm}=E_{2,2n}^{0,\pm}$.
Asymptotics \er{TM-4} follow.$\BBox$

Now we prove  Remark  to Theorem \ref{Tk}, \ref{4s}.

\begin{lemma} \lb{eq}
Let $k\ne{N\/2}$. Then for each $n:\vk_n\not\in R_k$ the following
relations hold true:
\[
\lb{eq1}
v_k(\l_{1,2n-1}^{0,\pm})\ge 0\qq\Leftrightarrow\qq
F_{k,1}(r_{k,n}^\pm)=F_{k,2}(r_{k,n}^\pm)\le -1
,
\]
\[
\lb{eq2}
v_k(\l_{1,2n-1}^{0,\pm})<0\ \ \Leftrightarrow\ \
F_{k,1}(r_{k,n}^\pm)=F_{k,2}(r_{k,n}^\pm)\in(-1,-{1\/2}].
\]
Moreover, if $q$ is even, i.e. $q(1-t)=q(t)$,
then $F_{k,1}(r_{k,n}^\pm)=F_{k,2}(r_{k,n}^\pm)\in(-1,-{1\/2}]$.
\end{lemma}

\no {\bf Proof.}
Let $r=r_{k,n}^\pm$. Recall $\r_k(r)=0$.
Identities \er{DeLk} yield
$(9F^2(r)-s_k^2)c_k^2=-s_k^2F_-^2(r)\le 0$, then
$9F^2(r)\le s_k^2$.
Moreover,
$$
F_{k,1}(r)=F_{k,2}(r)=\x_{k}(r)=
{9F^2(r)-F_-^2(r)-1\/2}-s_k^2\le
-{1\/2}.
$$
Relations \er{rg1}, \er{rg2} give that
if $v_k(\l_{1,p}^{0,\pm})> 0$, then
$F_{k,\n}(r)<-1$,
and if $v_k(\l_{1,p}^{0,\pm})<0$, then
$F_{k,\n}(r)>-1,\n=1,2$.

Conversely, let $F_{k,\n}(r)<-1$. Then \er{g1} yield $v_k(r)>0$.
Identities \er{sFc} give $v_k(r)>0$ on $\vk_{k,n}^\pm$, then
$v_k(\l_{1,p}^{0,\pm})\ge 0$.
Let $F_{k,\n}(r)>-1$. Identity \er{s1a} show that $9F^2(r)<h_1(r)$.
Then \er{g1} yield $v_k(r)<0$.
Identities \er{sFc} give $v_k(r)<0$ on $\vk_{k,n}^\pm$, then
$v_k(\l_{1,p}^{0,\pm})<0$.
Relations \er{eq1}, \er{eq2} are proved.

If $q\in L^2_{even}(0,1)$, then $F_-=0$ (see \cite{MW})
and $v_k<0,k\ne{N\/2}$.
Then \er{eq2} gives $F_{k,1}(r)=F_{k,2}(r)=
-{s_k^2+1\/2}\in(-1,-{1\/2}]$.
$\BBox$

\begin{proposition}\lb{pmg}
Let $k\not\in\{0,{N\/2}\}$, $q=q_\ve={1\/\ve}\d(t-{1\/2}-c_k \ve-\ve^2),\ve\ne 0$
and let $n_0>1$.
Then there exists $\ve_1>0$ such that for all
$\ve\in(-\ve_1,\ve_1)\sm\{0\},1\le n\le n_0$
the following relations hold true:
\[
\lb{vkd}
v_k(\l_{1,2n-1}^{0,\pm},q_\ve)>0,\qq
E_{1,2n-1}^{k,\pm}(q_\ve)=\l_{1,2n-1}^{0,\pm}(q_\ve),
\]
\[
\lb{uld}
u_\ell(E_{2,2n-1}^{0,\pm},q_\ve)>0,\ \
S_{1,n}^\ell(q_\ve)\cap S_{2,n}^\ell(q_\ve)\ne\es,\ \
\text{all}\qq \ 0\le\ell<{N\/2}-k.
\]
\end{proposition}

\no {\bf Proof}. If $q_v={1\/v}\d(t-a),v\ne 0,
a\in (0,1)$, then we have (see, for example \cite{BBKL})
$
F_-(\l,q_\ve)={\sin z(2a-1)\/2z\ve}, z=\sqrt\l.
$
Let $a={1\/2}+{c_k}\ve+{\ve^2},k\not\in\{0,{N\/2}\}$. Then
\[
\lb{P4b}
F_-(\l,q_\ve)={\sin 2z({c_k \ve}+{\ve^2})\/2z\ve}={c_k}+{\ve}+O(\ve^{2})\qq
\text{as}\qq |\ve|\to 0,
\]
uniformly on $|z|\le \pi n_0$.
Using this asymptotics we deduce that there exists $\ve_1>0$ such that
if $|\ve|<\ve_1$, then $|F_-(\l)|>c_k^2$ and $v_k(\l)>0$
for all $0\le\l<(\pi n_0)^2$. Using \er{s} we  obtain \er{vkd}.
Moreover, $|F_-(\l)|>s_\ell^2$ for all $\ell<{N\/2}-k$ and
$0\le\l<(\pi n_0)^2$. Thus, we obtain $u_\ell(\l)>0$ for such $\ell,\l$.
Then \er{mg} gives \er{uld}.
$\BBox$

 \no {\bf Acknowledgments.} 
Evgeny Korotyaev was partly supported by DFG project BR691/23-1.
The various parts of this paper were written at  ESI, Vienna , E. Korotyaev is grateful to the Institute for the hospitality.
A. Badanin is grateful to the Mathematical Institute of Humboldt Univ. for the hospitality.

\end{document}